\documentstyle[12pt]{article}
\begin{document}

\renewcommand{\theequation}{\thesection.\arabic{equation}}
%\footskip 17pt \addtolength{\leftmargin}{-0,5cm}
%\addtolength{\topmargin}{-0.5cm} \setlength{\textheight}{24cm}
%\setlength{\textwidth}{20,5cm} \baselineskip 10pt
%\setlength{\topmargin}{-30pt}

%\addtolength{\leftmargin}{-1,0cm} \addtolength{\topmargin}{-0.5cm}
\setlength{\textheight}{24cm} \setlength{\textwidth}{16,5cm}

\baselineskip 10pt \setlength{\topmargin}{-30pt}
\date{}
\begin{center}
{\bf \large Generalized solutions of the Cauchy problem for \\
 the Navier-Stokes system and diffusion
processes}\end{center}

 \vskip3mm
\begin{center}
S.~Albeverio \\
\small Institut f\"{u}r Angewandte Mathematik, Universit\"{a}t
Bonn,
 \\
\small Wegelerstr. 6, D-53115 Bonn, \ Germany\ SFB 611, \ Bonn, \ BiBoS, \\
Bielefeld - Bonn \\
\small \ CERFIM, Locarno and USI (Switzerland) \\
\end{center}
\begin{center}
Ya. Belopolskaya, \\
St.Petersburg State University
 for Architecture and Civil \\Engineering,
  2-ja Krasnoarmejskaja 4,\\ 190005, St.Petersburg, Russia
\end{center}
 %\end of header
\abstract{We reduce the construction of a weak solution
 of the Cauchy problem for the Navier-Stokes system on $R^3$ to the construction
 of a solution to a stochastic problem. Namely, we construct  diffusion
 processes which allow us to obtain a probabilistic representation of a weak 
 (in distributional sense)  solution to the Cauchy problem for the Navier- Stokes system on a small time interval. Strong solutions on a small time interval are constructed as well}

\vspace {2mm}

\vspace {2mm} {\bf Key words: } Stochastic  flows, diffusion
process, nonlinear para-\\bolic equations, Cauchy problem.
\section*{Introduction}
The main purpose of this article is to construct both  strong and
     weak solutions (in  certain  functional classes) of the Cauchy
problem for the Navier-Stokes (N-S) system in $R^3$. To this end we
consider a stochastic  problem and show that the solution of the
Cauchy problem for the Navier-Stokes system can be constructed via
the solution of this stochastic problem.

 The approach we develop  in this article is  based on the
 theory of stochastic equations associated with nonlinear parabolic
 equations started  by McKean \cite{MK} and Freidlin \cite{Fr},\cite{Fr1} and generalized by
 Belopolskaya and Dalecky  \cite{BD1}, \cite{BD2} on one hand and on the theory of
 stochastic flows due to Kunita \cite{Ku} on the other hand. In our previous
 paper \cite{AB} we have constructed a stochastic process that allows us to
 prove the existence and uniqueness of a local in time classical
 ($C^2$-smooth in the spatial variable) solution
 of the Cauchy problem for the Navier -Stokes
 system. In the present paper we construct a process which allows us to
 obtain construction of solutions of the both weak and strong Cauchy problem for
 this system. Later we plan to apply a similar approach for
 the Navier-Stokes equation for
 compressible fluids extending the results from \cite{AB2},
 \cite{AB3}.

A close but different approach is  the Euler-Lagrange approach to
incompressible fluids which was developed by Constantin \cite{Con1} and
Constantin and Iyer \cite{CI1}. Shortly, the main differences in these
approaches are the following: we use  a probabilistic
representation for the Euler pressure instead of  the Leray
projection and obtain different formulas for the stochastic
representation of the velocity field.   We discuss these
differences  with more details in the last section of the present work.

Thestructure of the present  article is  as follows. In the first section we give
some preliminary information concerning  different analytical
approaches to the  Navier-Stokes  system. Here we recall some
common ways to eliminate the pressure and to obtain a closed
equation for the velocity.

The classical approaches here are based on  the so called Leray
(Leray-Hodge)-projection that  is  a projection of the space of
square integrable vector fields to the space of divergence free
square integrable vector fields. Applying such a projection to the
velocity equation one can eliminate the pressure $p$ and get the
closed equation for the velocity $u$. This operator is used both
in numerous analytical papers (see \cite{LR} for references) and
in papers  where the  N-S system is studied from the probabilistic
point of view \cite{LS},\cite{Os}, \cite{CI1}. Finally the
pressure is reconstructed from the Poisson equation.

One more possibility to eliminate the pressure appears when one
considers the equation for the vorticity of the velocity field $u$
and  uses the Biot-Sawart law to obtain a closed system. From the
probabilistic point of view this approach was investigated in
\cite{BFR}.

In our previous paper \cite{AB} we do  not use the Leray
projection but instead we start with consideration of a system
consisting of the original velocity equation and the Poisson
equation for the pressure and construct their probabilistic
counterpart. The probabilistic counterpart of the N-S system was
presented in the form of a system of stochastic equations.
Furthermore we prove the existence and uniqueness  of a solution
to this stochastic system and show that in this way
 we  construct a unique classical (strong) solution of the Cauchy problem for the N-S system
  defined on a small time interval
 depending on the Cauchy data.

In the present paper we also reduce the N-S system to the system
of equations consisting of the original velocity equation and  the
Poisson equation for the pressure but then an associated stochastic
problem considered here allows to construct a generalized
(distributional) solution to the Cauchy problem for the N-S
system. The associated stochastic problem is studied in section 5.
In sections 1-4 we expose auxiliary results used in section 5.
Namely, in section 1 we give analytical preliminaries and recall
the notions of strong, weak and mild solutions to the Cauchy
problem for the Navier-Stokes system. More detail can be found the
recent book by Lemarie-Rieusset \cite{LR}. In section 2 we give a
short review of probabilistic approaches to the investigation of
the Navier-Stokes system \cite{AB}, \cite{LS} -\cite{BFR}. In
section 3 we study a probabilistic representation of the solution
to the Poisson equation, while in section 4 we recall some
principal fact of the Kunita theory of stochastic flows and apply
the results from \cite{Ku1}, \cite{Ku2} to construct a solution of
the Cauchy problem for a  nonlinear parabolic equation (see also
\cite{BW}). Finally all these preliminary results are used to
construct the probabilistic counterpart of the Navier-Stokes
system, prove that there exists a unique local solution to the
corresponding stochastic system and apply the results to construct
both the strong and weak (and simultaneously mild) solutions to
the Cauchy problem for the Navier-Stokes system.

\section{Preliminaries}

As it was mentioned in the introduction the main purpose of this
article is to construct both strong and
     weak solutions (in  certain  functional classes) of the Cauchy
problem for the Navier-Stokes  system  via diffusion processes.

Consider the Cauchy problem for the Navier-Stokes system
\begin{equation}\label{1.1}
 \frac{\partial u}{\partial t}+(u,\nabla )u=\nu\Delta u-\nabla p ,\quad u(0,x)=u_0(x),\quad x\in R^3,
\end{equation}
\begin{equation}\label{1.2}
 div\,u=0.\end{equation}
 Here $u(t,x)\in R^3,x\in R^3,t\in [0,\infty
)$ is the velocity of the fluid at the position $x$ at time $t$ and
 $\nu >0$ is the viscosity coefficient and  $p(t,x)$ is a scalar
field called the pressure which appears in the equation to enforce
the incompressibility condition (\ref{1.2}). Later we set
$\nu=\frac{\sigma ^2}2$ for reasons to be explained below.

By eliminating the  pressure from (\ref{1.1}),(\ref{1.2}) one  gets
a nonlinear pseudo-differential equation which is to be solved.
There exist different ways to do it and we consider now some of
them.

 Given a vector field $f$
let $\textbf{P} f$ be given by
\begin{equation}\label{1.4}
\textbf{P}f=f-\nabla \Delta^{-1} \nabla \cdot f.\end{equation}
Here and below we denote by $u\cdot v$ the inner product of
vectors $u$ and $v$ valued in $R^3$.

 The
map $\textbf{P}$  called  the  Leray projection is  a projection
of the space ${\bf L}^2(R^3)\equiv L^2 (R^3)^3$ of square
integrable vector fields to the space of divergence free vector
fields and we discuss its properties below. A quite direct
definition of $\textbf{P}$ is connected with the Riesz
transformation $R_j$. Recall that
$R_k=\frac{\nabla_k}{\sqrt{-\Delta}}$ which means that  for
$f\in{\bf L}^2$ we have ${\cal F}(R_jf)=\frac{i\xi_j}{|\xi|}\hat
f(\xi)$ where ${\cal F}(f)=\hat f$ is the Fourier transform of
$f$. Then $\textbf{P}$ is defined on  ${\bf L}^2(R^3)$ as
$\textbf{P}=Id+R\otimes R$ or
$$(\textbf{P}f)_j=f_j+\sum_{k=1}^3R_jR_kf_k.$$
Since $R_kR_j$ is a Calderon-Zygmund operator, $\textbf{P}f$ may
be defined on many Banach spaces.

Set
\begin{equation}\label{1.3}
\gamma(t,x)=\sum_{k,j=1}^3\nabla_ku_j\nabla_ju_k ={\rm Tr}[\nabla
u]^2\end{equation} and note that $\gamma$ can be presented as well
in the form
$$\gamma=\nabla \cdot\nabla \cdot u\otimes
u=\sum_{j,k}\nabla_k\nabla_j(u_ku_j).$$ .

 By computing
the divergence of both parts of (\ref{1.1}) and taking into
account (\ref{1.2}) we derive the equation
\begin{equation}\label{1.5}
-\Delta p(t,x)=\gamma (t,x)
\end{equation}
 thus arriving at the Poisson equation.
 The formal solution of  the Poisson equation is given by
\begin{equation}\label{1.6}
p=\Delta^{-1} \gamma = \Delta^{-1} \nabla \cdot\nabla \cdot
u\otimes u\end{equation}
since ${\rm div}u=0$
 and finally we present $\nabla p$ in the
form
$$\nabla p= \nabla \Delta^{-1} \nabla \cdot\nabla \cdot u\otimes
u.$$ Substituting this expression for $\nabla p$ into (\ref{1.1})
we obtain the following Cauchy problem
\begin{equation}\label{1.7}
\frac{\partial u}{\partial t}=\nu\Delta u-{\bf P}\nabla\cdot
(u\otimes u), \quad u(0)=u_0.\end{equation}

There are a number of ways to define a notion of a solution for
the Cauchy problem (\ref{1.7}).  We will
 appeal mainly to the Leray weak solution \cite{Ler} or to the Kato mild
solution \cite{Kato}.

\subsection{Leray and Kato approaches to the solutions of the Navier-Stokes equations}

Let ${\cal D}={\cal D}(R^3)=C_c^\infty$ be the space of all
infinitely differentiable functions on $R^3$  with compact support
equipped with the Schwartz topology. Let ${\cal D}'$ be the
topological dual of $\cal D$    and denote by $\langle
\phi,\psi\rangle=\int_{R^3}\phi(x)\psi(x)dx$ the natural coupling
between $ \phi\in{\cal D}$ and $\psi\in{\cal D}'$. If it will not
lead to misunderstandings we will use the same notation for vector
fields $u$ and $v$ as well, that is
$$\langle
h,u\rangle=\int_{R^3}\sum_{k=1}^3h_k(x)u_k(x)dx.$$

We recall that a weak solution of the N-S system on $[0,T]\times R^3$
is a distribution vector field $u(t,x)$ in $({\cal D}'((0,T)\times
R^3))^3$ where $u$ is locally square integrable on $(0,T)\times
R^3$,  div $u=0$ and there exists $p\in {\cal D}'((0,T)\times
R^3)$ such that
\begin{equation}\label{1.0}
 \frac{\partial u}{\partial t}=\nu\Delta u- \nabla\cdot(u\otimes u) - \nabla p ,\quad
 \lim_{t\to 0} u(t)=u_0
\end{equation}
holds.

 The Leray
solution to the N-S equations is constructed through a limiting
procedure from the solutions to the mollified  N-S equations
\begin{equation}\label{1.8}
\cases{\frac{\partial u}{\partial t}= \nu\Delta u-\nabla
\cdot((u*q_\varepsilon)\otimes u)-\nabla p,\cr \nabla\cdot u=0,\cr
\lim_{t\to 0}u(t)=u_0.}
\end{equation}
Namely it is proved that  there exists a function
$$u_\varepsilon\in L^\infty((0,\infty),{\bf L}^2)\cap
L^2((0,T),(\dot {\bf H}^1))$$ such that (at least for a
subsequence $u_{\varepsilon_k}$) strongly converging  in
$(L_{loc}^2((0,T)\times R^3))^3$ to $u$ which satisfies
(\ref{1.8}).

Here $\dot{\bf H}^1$ is the homogenous Sobolev space $\dot{\bf
H}^1=\{f\in {\bf S}_0':\nabla f\in{\bf L}^2\}$ with norm
$\|f\|_{{\bf H}_1}=\|\nabla f\|_{{\bf L}_2}$.

On the other hand to construct the Kato solution means to
construct a solution $u$ to the following integral equation
\begin{equation}\label{1.9}
u(t)=e^{t\Delta}u_0-\int_0^te^{(t-s)\Delta}{\textbf
P}\nabla\cdot(u\otimes u)(s)ds.
\end{equation}

Note that instead of looking for $u(t,x)$ and $p(t,x)$ one can
prefer to look for their Fourier images $\hat
u(t,\lambda)=(2\pi)^{-\frac 32} \int_{R^3}e^{-i\lambda\cdot
x}u(t,x)dx$.

The Leray and Kato approaches stated in terms of the Fourier
transformations of the Navier-Stokes system can be described as
follows.

Applying the Fourier transformation to the relation (\ref{1.7})
written in the form
$$\langle h, u(t)\rangle=\langle h, u(0)\rangle+\int_0^t\langle h,\Delta u(s)\rangle-
\int_0^t\langle h,\nabla\cdot(u\otimes u)(s)\rangle$$ we derive
the relation
\begin{equation}\label{1.10}
\langle \hat h, \hat u\rangle= \langle \hat h, \hat u_0\rangle
-\int_0^t\langle \hat h, |\lambda|^2\hat u(s)\rangle
ds-\end{equation}
$$ \frac{i}{(2\pi)^{\frac
32}}\int_0^t\int_{R^3}\int_{R^3} \sum_{k,l=1}^3\lambda^k \hat
h^l(\lambda'), \hat u_l (s,\lambda)
u_k(s,\lambda-\lambda')d\lambda d\lambda'ds.$$ Here $\hat u$
corresponds to the Fourier transformation of $u$.

On the other hand if we are interested in the Kato mild solution
of the N-S system then we may apply the Fourier transformation to
(\ref{1.9}) and derive the following equation
\begin{equation}\label{1.11}
\chi(t,\lambda)=exp\{-\nu|\lambda|^2t\}\chi(0,\lambda)+\end{equation}
$$\int_0^t\nu|\lambda|^2
e^{-\nu|\lambda|^2(t-s)}\left[\frac 12
(\chi(s)\circ\chi(s))\right](\lambda)ds$$ for the function
$$\chi(t,\lambda)=\frac 2\nu\left(\frac \pi 2\right)^{\frac
32}|\lambda|^2\hat u(t,\lambda).$$
Here
\begin{equation}\label{1.12}
\chi_1\circ\chi_2(\lambda)=-\frac{i}{\pi^3}\int_{R^3}
(\chi_1(\lambda_1)\cdot e_\lambda)
\Pi(\lambda)\chi_2(\lambda-\lambda')\frac{|\lambda|d\lambda'}{|\lambda'|^2|\lambda-\lambda'|^2},
\end{equation}
$e_\lambda=\frac{\lambda}{|\lambda|}$ and
\begin{equation}\label{1.13}
\Pi(\lambda)\chi=\chi-e_\lambda(\chi\cdot e_\lambda),
\end{equation}

Coming back to (\ref{1.7}) we note  that the Leray projection
allows to  eliminate  the pressure $p(t,x)$ from the Navier-Stokes system,
to construct $u$ and finally to look for $p$ defined by the
solution of the auxiliary Poisson equation.

Another way to eliminate $p(t,x)$ from the system
(\ref{1.1}),(\ref{1.2}) is to  consider the function $v(t,x)={\rm
curl}\, u(t,x)$ called the vorticity. Since ${\rm curl}\nabla
p(t,x)=0$ one can derive a closed system for $u$ and $v$.  Namely
for $u$ and $v$ we arrive at the system consisting of the equation
\begin{equation}\label{1.14}\frac{\partial v}{\partial t}+(u\cdot\nabla)v=\nu \Delta
v+(v\cdot\nabla)u ,\end{equation} and the so called Biot-Savart
law having the form
\begin{equation}\label{1.15}u(t,x)= \frac 1{4\pi}\int_{R^3}\frac{(x-y)\times
v(y)}{|x-y|^3}dy.
\end{equation}
Here the cross-product $ u\times v$ is given by
$$u\times v=det\pmatrix{e_1&e_2&e_3\cr u_1&u_2&u_3\cr v_1&v_2&v_3\cr}=$$
$$(u_2v_3-u_3v_2)e_1+(u_3v_1-u_1v_3)e_2+(u_1v_2-u_2v_1)e_3,$$
where $(e_1,e_2,e_3)$ is the orthonormal basis in $R^3$.

Note that the term $(v\cdot\nabla)u$ can be  written as $(\nabla
u)v$ or even as ${\cal D}_uv,$ where ${\cal D}_u$ is the
deformation tensor defined as  the symmetric part of $\nabla u$
$${\cal D}_u=\frac 12(\nabla\, u+\nabla\, u^T),$$
since by direct computation we see that
$$(\nabla \,u)v-{\cal D}_uv=\frac 12(\nabla\, u+\nabla\,
u^T)v=0.$$

To be able to present the precise statements concerning the
existence and uniqueness of solutions to the N-S equations we have
to introduce a number of functional spaces to be used in the sequel.
\subsection{Functional spaces}
We describe here functional spaces which will be used in the
sequel.

Let ${\cal D}={\cal D}(R^3)$ be the space of all infinitely differentiable functions on $R^3$ with compact supports equipped with the Schwartz topology. Let ${\cal D}'$ be the topological dual to ${\cal D}$. The elements of ${\cal D}'$  are called Schwartz distributions.

The space of $R^3$-valued vector fields  $h$   with
components $ h_k\in {\cal D}$ shall be  denoted  by  ${\bf D}(R^3)$ and
  ${\bf D}'$  shall denote  the space  dual to ${\bf D}(R^3)$.

Let  $L^q(R^3)$ denote the Banach space of functions $f$ which are  absolutely
integrable taken to the $q$-th power with the norm
$\|f\|_q=(\int_{R^3}|f(x)|^qdx)^{\frac 1q}$;

Let $Z$ denote the set of all integers, and suppose that $k\in Z$
is  positive and $1<q<\infty$. Denote by
 $W^{k,q}=W^{k,q}(R^3)$ the set of
 all real functions  $h$ defined on $R^3$ such that $h$ and all
 its distributional derivatives $\nabla^\alpha$ of order
 $|\alpha|=\sum\alpha_j\le k$  belong to $L^q(R^3)$.  It is a
 Banach space with norm
\begin{equation}\label{1.16}
\|h\|_{k,p}=(\sum_{|\alpha|\le k}\int_{R^3}|D^\alpha
h(x)|^qdx)^{\frac 1q}.
\end{equation}
We denote the dual space of $W^{k,q}$ by $W^{-k,m}$ where $\frac
1m+\frac 1q=1$. Elements of $W^{-k,q}$ can be identified with
Schwartz distributions. The space $W^{-k,q}$ is also a Banach
space with norm
$$\|\phi\|_{-k,q}=\sup_{\|h\|_{k,q}\le 1}|\langle
\phi, h\rangle|, $$ where
$$\langle
\phi, h\rangle=\int_{R^3}\phi(x)h(x)dx.$$ The spaces $W^{k,p}$ for
$k\in Z$ and $p>1$ are called Sobolev spaces. If $p=2$ we use the
notation $H^k$ for the  Hilbert spaces $W^{k,2}$ .
In a natural way one can define the spaces $\bf{W}^{k,q}$, ${\bf
H}^k$ of vector  fields with components in $W^{k,p}$, and $H^k$
and so on.

Set
$${\cal V}=\{v\in {\bf D}: div v=0\}$$
 and let
\begin{equation}\label{1.17}
{\bf H}=\{\mbox{closure of } \,{\cal V}\quad \mbox{in }\, {\bf
L}^2(R^3)\}, \quad {\bf V}=\{\mbox{closure of }\, {\cal V}\quad
\mbox{in }\, {\bf H}^1\}.
\end{equation}

Let   $C_b^k(R^3, R^3)$ denote the space of k-times differentiable
fields with the norm
$$\|g\|_{C^k_b}=\sum_{|\beta|\le k}\|D^\beta g\|_\infty$$
and let $C_b^{k,\alpha}(R^3, R^3)$ be the space of vector fields
whose k-th derivatives are H\"older continuous with exponent
$\alpha,\,$ $0<\alpha<1$ with the norm
$$\|g\|_{C_b^{k,\alpha}}=\|g\|_{C_b^k}+[g]_{k+\alpha}$$
where
$$[g]_{k+\alpha}=\sum_{|\beta|=k}\sup_{x,y\in R^3}\frac{|D^\beta g(x)
-D^\beta g(y)|}{|x-y|^\alpha}.$$

We denote by Lip$(R^3)$ the space of bounded Lipschitz continuous
functions with the norm
$$\|g\|_{Lip}=\sup_{x,y\in R^3}\frac{|g(x)
-g(y)|}{|x-y|}.$$

Spaces of integrable functions on the whole $R^3$ appear to be
not  satisfactory to construct a solution to the N-S equations and
one has to consider spaces of locally integrable functions.

Let $f: R^3\to R^1$ be a Lebesgue
measurable function. A set of functions $\{f:
\int_{K}|f(x)|^pdx<\infty\}$  for all compact subsets $K$ in $R^3$
is denoted by $ L^p_{loc}$ and called a space of locally
integrable functions.  Note that $L^1(R^3)\subset L^1_{loc}(R^3)$ .  Although $L^p_{loc}(R^3)$ are not normed spaces
they are readily topologized. Namely a sequence $\{u_n\}$
converges to $u$ in  $L^p_{loc}(R^3)$ if $\{u_n\}\to
u$  in $L^p(K)$ for each open $K\subset G$  having  compact closure in $R^3$.
 Local spaces $W^{k,p}_{loc}(R^3)$ can be defined to consist of
 functions belonging to $W^{k,p}(K)$ for all compact $K\subset R^3$.

A local space $W^{k,p}_{loc}(G)$  is defined as a space of functions
belonging to $W^{k,p}(G')$  for all $G'\subset G$ with compact
closure in $G$.  A  function $f\in W^{k,p}_{loc}(G)$  with compact
support will in fact belong to $W^{k,p}_0(G).$ Also  functions in
$W^{1,p}(G)$ which vanish continuously on the boundary $\partial G$
will  belong to $W^{1,p}_0(G)$  since they can  be approximated by
functions with compact support.

 In the whole space $R^3$ and with $p$, $q$ satisfying $1\le q\le
p<\infty$ denote  by ${\cal M}^p_q$ a nonhomogenous Morrey space  and by
$ M^p_q$ a homogenous Morrey space with norms given respectively by

\begin{equation}\label{1.19}
{\cal M}_q^p=\left\{f\in L^q_{loc}:\|f\|_{{\cal
M}^p_q}=\sup_{x_0\in R^3}\sup_{0<R}R^{\frac 3p-\frac
3q}\|f\|_{L^q(B(x_0,R))}<\infty\right\},
\end{equation}
and

\begin{equation}\label{1.20}
{ M}_q^p=\left\{f\in L^q_{loc}:\|f\|_{{\cal M}^p_q}=\sup_{x_0\in
R^3}\sup_{0<R\le 1}R^{\frac 3p-\frac
3q}\|f\|_{L^q(B(x_0,R))}<\infty\right\}
\end{equation}
where $B(x_0,R)$ is a closed ball of $R^3$ with center at $x_0$
and radius $R$.

Respectively the integrable function is said to belong to $M^q(G)$ if there exists a constant C such that
\begin{equation}\label{1.02}
\int_{G\cap B_R}|f(x)|dx\le CR^{3(1-\frac 1q)}\end{equation}
for all balls $B_R$. The norm in $M^q(G)$  is defined as the minimum of the constants $C$
satisfying (\ref{1.02})

A distribution $u$ on $(0,T)\times R^3$ is said to be uniformly
locally square integrable if for all $\varphi\in {\cal
D}((0,T)\times R^3) $
$$\sup_{x_0\in R^3} \int_0^T\int_{R^3}|\|\varphi(t,
x-x_0)u(t,x)\|^2dxdt<\infty.$$

Equivalently $ u$ is uniformly locally square integrable if and
only if for all $t_0<t_1\in (0,T)$ the function
$U_{t_0,t_1}(x)=(\int_{t_0}^{t_1}\|u(t,x)\|^2dt)^{\frac 12}$
belongs to the Morrey space $L^2_{uloc}$ . In this case we write
$$u\in \cap_{0<t_0<t_1<T}({\bf L}^2_{uloc}{\bf L}^2_t((t_0,t_1)\times
R^3)).$$

  For $1\le p\le \infty $ the Morrey space of uniformly
locally integrable functions on $R^3$ is the Banach space
$L^p_{uloc}$ of Lebesgue measurable functions $f$ on $R^3$ such that
the norm $\|f\|_{p, uloc}$ is finite, where
$$\|f\|_{p, uloc}=\sup_{x_0\in
R^3}(\int_{\|x-x_0\|<1}|f(x)|^pdx)^{\frac 1p}.$$

For $t_0<t_1,\, 1\le p,q\le \infty$ the space $L^p_{uloc,
x}L^q_t((t_0,t_1)\times R^3))$ is the Banach space of  Lebesgue
measurable functions $f$ on $(t_0, t_1)\times R^3$ such that the
norm
$$\sup_{x_0\in
R^3}(\int_{\|x-x_0\|<1}(\int_{t_0}^{t_1}|f(t,x)|^qdt)^{\frac
pq}dx)^{\frac 1p}$$ is finite.

A $C^\infty$ function   $f$ on $R^3$ is called  rapidly decreasing
if  $$\lim_{x\to\infty}|D^\alpha f(x)|(1+\|x\|)^n=0$$ holds for
any multi-index $\alpha$ and any positive integer $n$.  Let ${\cal
S}={\cal S}(R^3)$ be the space of rapidly decreasing $C^\infty-$
functions equipped with the Schwartz topology and ${\cal S}'$ be
the topological dual of $\cal S$. Since $\cal S$ includes $\cal
D$, ${\cal S}'$ is a subset of ${\cal D}'$. The elements of ${\cal
S}'$ are  called tempered distributions.

\subsection{Weak, strong  and mild solutions of the Navier-Stokes
system}

 Now we are ready to give more precise definitions and
statements concerning the existence and uniqueness of solutions of
the N-S equations.

{\bf Definition 1.1.(Weak solutions)}{\em A weak solution of the
Navier-Stokes system on $(0,T)\times R^3$ is a distribution vector
field $u(t,x)$,  $u\in ({\bf {\cal D}}'((0,T)\times R^3)^d$ such
that

a) $u$ is locally square integrable on $(0,T)\times R^3$ ,

b) $\nabla\cdot u=0,$

c) there exists $p\in{\cal D}'((0,T)\times R^3)$ such that
$$\partial_t u=\Delta u-\nabla\cdot(u\otimes u)-\nabla p.$$}

The classical results concerning the existence of square
integrable weak solutions are due to Leray  \cite{Ler}.

{\bf Theorem 1.1.}   (Leray's theorem) {\em Let $u_0\in
(L^2(R^3))^3$ so that $\nabla\cdot u=0$. Then there exists a weak
solution $u\in L^\infty((0,\infty), (L^2)^d)\cap L^2((0,\infty),
(H^1)^3)$ for the Navier -Stokes equation on $ (0,\infty)\times
R^3$ so that $lim_{t\to 0}\|u(t)-u_0\|_2=0$. Moreover, the
solution $u$ satisfies the energy inequality
\begin{equation}\label{1.21}
\|u(t)\|^2+ 2 \int_0^t\int_{R^3}\|\nabla\otimes u\|^2dxds\le
\|u_0\|^2_2,
\end{equation}
where
$$\|u(t)\|^2_2=\sum_{k=1}^d\int_{R^3}|u_k(t,x)|^2dx,\quad
\nabla\otimes u=\sum_{k=1}^d\sum_{j=1}^d|\partial_ku_j|^2.$$}

{\bf Definition 1.2. (Mild solution)} {\em The Kato  mild solution
of (\ref{1.1}), (\ref{1.2}) is a solution of (\ref{1.7})
constructed as a fixed point of the transform
\begin{equation}\label{1.22}
v\mapsto e^{t\Delta} u_0(x)-\int_0^te^{(t-s)\Delta}{\bf
P}\nabla\cdot (v\otimes v)(\theta,x)d\theta  =e^{t\Delta}u_0 -
B(v,v).
\end{equation}. }

Note that  the right hand side of (\ref{1.7})
\begin{equation}\label{1.23}
e^{t\Delta} u_0(x)-\int_0^te^{(t-s)\Delta}{\bf P}\nabla\cdot
(u\otimes u)(\theta,x)d\theta  =\Phi(t,x,u)
\end{equation}
is  a nonlinear map in the corresponding space  and  the solution
$u$ is obtained by the iterative procedure
\begin{equation}\label{1.24}u^0=e^{t\Delta}u_0,
 \quad u^{n+1}=e^{t\Delta}u_0-B(u^n,u^n).\end{equation}

Hence to construct a mild solution to (\ref{1.1}), (\ref{1.2})
means to find a suitable functional space for which $\Phi(t,x,u)$
given by (\ref{1.20}) is a contraction.

To this end one has to find a subspace  ${\cal E}_T$ of
$L^2_{uloc, x}L^2_t((0,T)\times R^3)$ so that the bilinear
transformation $B(u,v)$ of the form (\ref{1.10}) is bounded as a
map ${\cal E}_T\times {\cal E}_T\to {\cal E}_T$. Then one may
consider the space ${\bf E}_T\subset {\cal S}'$ defined by $f\in
{\bf E}$ iff $ f\in {\cal S}'$ and $(e^{t\Delta}f) _{0<t<T}\in
{\cal E}_T$ and prove the following result.

{\bf Theorem 1.3}. The Picard contraction principle.

{\em Let ${\cal E}_T\subset L^2_{uloc,x}L^2_t([0,T)\times R^3)$ be
such that the bilinear map $B$ is bounded on  ${\cal E}_T$ Then:

(a) If $u\in {\cal E}_T$ is a weak solution for the Navier-Stokes
equation (\ref{1.1}) (\ref{1.2}) then the associated initial value
belongs to ${\bf E}_T$.

(b) There exists a positive constant  $C$ such that for all
$u_0\in {\bf E}_T$ satisfying $\nabla\cdot u=0$  and
$\|e^{t\Delta} u_0\|_{{\cal E}_T}<\infty$ there exists a weak solution
$u\in {\cal E}_T$ of  (\ref{1.1}) (\ref{1.2}) associated with the
initial value $u_0$
\begin{equation}\label{1.25}
u=e^{t\Delta}u_0-\int_0^te^{(t-s)\Delta}{\bf
P}\nabla\cdot(u\otimes u)ds.
\end{equation}}

The classical results assert  that  for sufficiently smooth
initial data for example for $u_0$ in the Sobolev space $H^k$ ,
$k>\frac d2+1,$ there exists a short time strong unique solution to
({\ref{1.1}), ({\ref{1.2}).  On the other hand Leray  proved the
existence of a global weak solution of finite energy, i.e. $u\in L^2$
called the Leray-Hopf solution. Although the uniqueness and full
regularity of this solution are still an open problem nevertheless one knows that
if a strong solution exists then a weak solution coincides with
it.

\section{Probabilistic  approaches to the solution of the N-S equations  }
\setcounter{equation}{0}

 Along with above
functional analytical approaches recently a number of
probabilistic approaches to the problems of hydrodynamics was
developed (see \cite{BFR}-\cite{Os}, \cite{AB}). In this section
we give a short survey of several different probabilistic
approaches.

Let $(\Omega,{\cal F},P)$ be a complete probability
space, $w(t), B(t)$ be a couple of independent Wiener processes
valued in $R^3$.

Assume that  $(u(t,x), p(t,x))$ is a unique strong solution to
(\ref{1.1}), (\ref{1.2}) or (to be more precise) to
(\ref{1.1}), (\ref{1.5}) and set $\nu=\frac{\sigma^2} 2$ to
simplify notations in  stochastic equations. Since in this case
$u$ is $C^2$-smooth  one can check that the stochastic equation
\begin{equation}\label{2.1}
d\xi(\tau)=-u(t-\tau,\xi(\tau))d\tau+\sigma dw(\tau),\quad
\xi(0)=x,
\end{equation}
has a unique solution and the function
\begin{equation}\label{2.2}
v(t,x)=E[u_0(\xi(t))-\int_0^t\nabla p(t-\tau,\xi(\tau))d\tau]
\end{equation}
satisfies (\ref{1.1}) and hence equals to $u$ by the uniqueness
of the strong solution to (\ref{1.1}). The relation
\begin{equation}\label{2.3}
-2p(t,x)=\int_0^\infty E\gamma(t,x+B(\theta))d\theta
\end{equation}
with $\gamma $ given by (\ref{1.3}) allows us to verify that ${\rm
div}u=0$. Thus (\ref{2.2}),(\ref{2.3}) give the
probabilistic representation of the solution to the N-S system.

  If $u(t,x)$ is a $C^2$-smooth solution to
(\ref{2.1})-(\ref{2.3}), then Ito's formula yields that for
$v(\theta, x)=u(t-\theta, x)$
$$ v(t,\xi (t))=v(\theta,x)+\int^t_\theta [\frac{\partial v}{\partial \tau
}-(v,\nabla )v+\frac{\sigma ^2}2\Delta v](\tau ,\xi (\tau ))d\tau
+$$$$\int^t_\theta \sigma \nabla v(\tau,\xi(\tau))dw(\tau ).
$$
Then (\ref{2.2}) and the relation $u_\tau^{\prime }(t-\tau
,x)=-v_\tau ^{\prime }(\tau ,x)$ yield
$$
Eu(0, \xi(t))= v(t,x)-E[\int^t_0 [\frac{\partial u}{\partial \tau
}+(u,\nabla )u-\frac{\sigma ^2}2\Delta u](t-\tau ,\xi_x (\tau
))+$$$$\nabla p(t-\tau ,\xi_x(\tau ))]d\tau]+E[\int^t_0\nabla
p(t-\tau ,\xi_x(\tau ))]d\tau].
$$
Finally it results from (\ref{2.2}) that
$$E[\int^t_0 [\frac{\partial u}{\partial \tau
}+(u,\nabla )u-\frac{\sigma ^2}2\Delta u+\nabla p](t-\tau
,\xi_x(\tau ))]d\tau]=0.$$
 Since the latter equality 
holds for all $t$ and $x$ we deduce that (\ref{1.1})also holds. The
relation $u(0,x)=u_0(x)$ immediately follows from (\ref{2.2}).

 The system (\ref{2.1})-(\ref{2.3}) is a closed system of
equations and  we can  try to look for  its solution. Then at a
 second step we will look for the connection between this solution and a solution of
  the N-S system.

This approach was realized in paper \cite{AB}. It appears that to
prove the existence of smooth solutions to (\ref{2.1})-(\ref{2.3})
 we have to consider the stochastic  representations
for $\nabla u$ and $\nabla p$ along with  this system.

 Using general results of diffusion
process theory and in particular the Bismut-Elworthy formula \cite{EL} we note
that heuristic differentiation of (\ref{2.1})- (\ref{2.3}) leads to
the relations
$$
\nabla _ku_i(t,x)=E[\nabla _ju_{0i}(\xi (t))\eta _{jk}(t)-$$
\begin{equation}\label{2.4}
\int_0^t\frac 1{\sigma \tau }(\nabla _ip(t-\tau ,\xi (\tau
))\int_0^\tau \eta _{kl}(\theta )dw_l(\theta ))d\tau ]
\end{equation}
and
\begin{equation}\label{2.5}
d\eta _{ik}=-\nabla _ju_i(t-\tau ,\xi (\tau ))\eta _{jk}(\tau )d\tau
,\quad \eta_{ik}(0)=\delta _{ik}.
\end{equation}

In addition  by Bismut-Elworthy's formula (integration by parts) we
can derive from (\ref{1.16}) the probabilistic representation for
$\nabla p(t,x)$
\begin{equation}\label{2.6}
 2\nabla p(t,x)=-\int_0^\infty \frac
1sE[\gamma (t,x+B(s))B(s)]ds.
\end{equation}

The main results in \cite{AB} can be stated in  the following way.

Let $V=(u, \nabla u)$, ${\cal V}_1=
\{V(t,x):\|V(t)\|_{L^r_{loc}}<\infty\},$ if $ \frac
53<r<2$ and
 ${\cal V}_2=
\{V(t,x):\|V(t)\|_{L^r_{loc}}<\infty\},$  if  $r>3$,
${\cal V} ={\cal V}_1\cap {\cal V}_2\cap C^{1+\alpha}$ and  let ${\cal
M}=C([0,T], {\cal V})$ denote the Banach space with the norm
$$\|V\|_{r,\alpha}= sup_{t\in [0,T]}[\|{ V}(t)\|_{\cal V} + [\nabla
u(t)]_\alpha].$$

{\bf Theorem 2.1.}({\cite{AB}){\em Assume that $V(0)=V_0\in {\cal V}$. Then there
exist a bounded interval $[0,T]$ depending on $V_0$  and a unique
solution $(\xi(t), u(t,x), p(t,x),$ $ \eta(t), \nabla u(t,x))$ to
the system (\ref{2.1})-(\ref{2.5}) belonging to ${\cal M}$ for each $\tau\in
[0,T]$.}

{\bf Theorem 2.2.}({\cite{AB}){\em Assume that the conditions of theorem 2.1 hold
and $u_0\in C^{2+\alpha}$. Then there exists  an  interval
$[0,T_1],$ $T_1\le T$, such that for all $t\in[0,T_1]$  there exists
a unique solution to (\ref{1.1}), (\ref{1.3}) in $\tilde{\cal
M}\subset{\cal M}$ where $\tilde{\cal M}={\cal M}\cap C^2$  and this solution is given by (\ref{2.2}),(\ref{2.3}).}

 A close approach based on a similar diffusion process was developed by Busnello,
Flandoli, Romito  \cite{BFR}, though their starting point  was the
system  that governs the vorticity $v={\rm curl}\, u$ and  velocity
$u$. The corresponding probabilistic counterpart  of
(\ref{1.14}),(\ref{1.15}) can be presented in the form of the
following stochastic system

\begin{equation}\label{2.7}
d\xi(\tau)=-u(t-\tau,\xi(\tau))d\tau+\sigma dw(t),\,\xi(s)=x,
\end{equation}
 and the following two relations
\begin{equation}\label{2.8}v(t,x)=E[U(t,s)
u_0(\xi(t)),\end{equation}
\begin{equation}\label{2.9}2u(t,x)=\int_0^\infty \frac 1 \theta E[v(t,x+ B(\theta))\times
B(\theta)]d\theta ,\end{equation} where $ U(t,s)=exp(\int_s^t\nabla
u(t-\tau,\xi(\tau)d\tau)$ is a solution to the linear equation
\begin{equation}\label{2.10} dU_s^{t,x}=\nabla u(t-s,\xi(s))U_s^{t,x}ds,\quad
U_0^{t,x}=Id.\end{equation}

The main results in \cite{BFR} are as follows.  Denote by
$${\cal U}^\alpha(T)=$$$$\{u\in C([0,T], C_b^1(R^3,
R^3))\cap L^\infty([0,T], C_b^{1,\alpha}(R^3,R^3))|{\rm
div}\,u=0\}
$$
the Banach space endowed with the norm
$$\|u\|_{{\cal U}^\alpha}=\sup\,{\rm ess}_{0\le t\le
T}\|u(t)\|_{C_b^{1,\alpha}}$$ and by
$${\cal V}^{\alpha,q}(T)=\{v\in C([0,T], C_b(R^3,
R^3))\cap L^\infty([0,T], C_b^\alpha(R^3,R^3))\}
$$
the Banach space endowed with the norm
$$\|v\|_{{\cal V}^{\alpha,q}}=\sup\,{\rm ess}_{0\le t\le
T}\|v(t)\|_{L^q\cap C_b^{\alpha}}.$$

 {\bf Theorem 2.3.}(\cite{BFR}) {\em Given
$p\in[1,\frac 32), \alpha\in [0,1]$ and $T>0$ let $\xi_0\in
C^\alpha_b(R^3,R^3)\cap L^p(R^3,R^3)$ and set
$$\varepsilon_0=\|v_0\|_{C^\alpha_b\cap L^p}.$$
Then there exists $\tau\in[0,T]$ depending only on
$\varepsilon_0$, such that there is a unique solution of
(\ref{2.7})-(\ref{2.10}). The diffusion process $\xi(t)$
plays a role of a  Lagrangian path, vector field  $u$ belongs to ${\cal
U}_\alpha,\,$ and vector field $ v$ belongs to $ {\cal V}^{\alpha, p}(\tau)$. In addition
the deformation matrix $U_s^{x,t}$ satisfies (\ref{2.7}).}

After these developments P.Constantin kindly attracted our attention to his papers \cite{Con1} \cite{CI1} 
 where the Lagrangian  approach was successfully applied to the investigation of
 the Navier-Stokes system. The presentation of  this approach
  and  the discussion  of their similarity and difference will be given in the last section of the present article.

A probabilistic representation of the solution to the Fourier
transformed Navier-Stokes equation (\ref{1.10}) was constructed by
Le Jan and Sznitman \cite{LS}.

To describe their approach recall the definition of the solution of the Fourier
representation  (FNS) of the Navier-Stokes system.

First for a solution $u(t,x)$ of (\ref{1.25}) and its Fourier
transform $\hat u(t,\lambda)$ one can introduce a
function $\chi(t,\lambda)$ defined on $[0,T]\times R^3$  such that
$$\chi_t(\lambda)=\frac 2\nu (\frac\pi 2)^{\frac 32}|\lambda|^2\hat
u_t(\lambda), \mbox { a.e. for } t\in[0,T],$$ and
$$ \chi_t(\lambda)\cdot \lambda=0,\quad \chi_t(-\lambda)=\bar{\chi}_t(\lambda).$$
In addition,  for Lebesgue a.e. $\lambda$,  $\chi_t(\lambda)$ solves the
equation
\begin{equation}\label{2.11}
\chi_t(\lambda)=exp(-\nu|\lambda|^2t)\chi_0(\lambda)+\int_0^t
\nu|\lambda|^2e^{-\nu|\lambda|^2(t-s)}\frac 12[
\chi_s\circ\chi_s(\lambda)]ds,\end{equation}
 where
$$\chi_1\circ\chi_2(\lambda)=-\frac i{\pi^3}\int (\chi_1(\lambda_1)\cdot
e_\lambda)\Pi(\lambda)\chi_2(\lambda-\lambda_1)\frac{
|\lambda|d\lambda_1}{|\lambda_1|^2|\lambda-\lambda_1|^2},$$ if the
initial function $\chi_0$ is a measurable function and for a.e.
$\lambda\in R^3\backslash\{0\}$
$$
\chi_0:R^3\backslash\{0\}\to C^3,\quad
\chi_0(\lambda)\cdot\lambda=0,\quad
\chi_0(-\lambda)=\overline{\chi_0(\lambda)}\, .
$$

Introducing the kernel $K$  from $R^3\backslash {0}$ to
$(R^3\backslash {0})^2$
$$\int h(\lambda_1,\lambda_2)K_\lambda(d\lambda_1,
d\lambda_2)=\frac{1}{\pi^3}\int
h(\lambda_1,\lambda-\lambda_1)\frac{|\lambda|d\lambda_1}{|\lambda_1|^2|\lambda-\lambda_1|^2}$$
for $h\ge 0$ measurable on $(R^3\backslash {0})^2$ one gets
$$\chi_1\circ\chi_2(\lambda)=-i\int(\chi_1(\lambda_1)\cdot
e_\lambda)\Pi(\lambda)\chi_2(\lambda_2)K_\lambda(d\lambda_1,d\lambda_2).$$

It turns out that $K$ is a Markovian kernel with some remarkable
features that allow to study  existence and uniqueness problems
for  (\ref{2.11})  with the help of a critical branching process on
$R^3\backslash \{0\}$ called the stochastic cascade.
Namely, LeJan and Sznitman have described a particle located in $\lambda$ such
that after an exponentially holding time with parameter
$\nu|\lambda|^2$ with equal probability the particle either dies or gives birth
to two descendants, distributed according  to $K_\lambda$. A
representation formula for the solution of (\ref{2.11}) is
constructed as the expectation of the result of a certain
operation performed along the branching tree generated by the
stochastic cascade.

 One more probabilistic model  was recently developed by  M.
Ossiander \cite{Os}.    A binary branching process with jumps that
corresponds to the formulation of solutions to N-S in physical
space was constructed in \cite{Os}.

Once again the N-S system is reformulated  incorporating
incompressibility via the Leray projection $\textbf{P}$ and then
the Duhamel principle is applied to derive
\begin{equation}\label{2.16}u=e^{\nu
t\Delta}u_0-\int_0^t e^{-\nu(t-s)\Delta}\textbf{P}\nabla\cdot(u\otimes u)(s)ds+\int_0^te^{-\nu(t-s)\Delta}{\cal
P}f(s)ds,\end{equation}
\begin{equation}\label{2.17}\nabla\cdot
u_0=0.\end{equation}

Let $K(y,t)=(2\pi t)^{-\frac 32}e^{-\frac{|y|^2}{2t}}$ be the
transition density of the Brownian motion $w(t)\in R^3$
$$({\bf P}_y)_{ij}=\delta_{ij}-(e_y)_i (e_y)_j$$
$$b_1(y;u,v)=(u\cdot e_y){\bf P}_yv+(v\cdot e_y){\bf P}_yu,$$
$$b_2(y;u,v)=b_1(y;u,v)+u\cdot(I-3e_ye_y^T)ve_y.$$
Then  (\ref{2.16}) can be rewritten in the form
\begin{equation}\label{2.18}
u(x,t)=\int_{R^3}u_0(x-y)K(y, 2\nu t)dy+$$
$$\int_0^t\int_{R^3}\{\frac{|z|}{4\nu s}K(z, 2\nu s)b_1(z;u(x-z,t-s),
u(x-z,t-s))+$$
$$(\frac 1{|z|}K(z, 2\nu s)-\frac
3{4\pi|z|^4}\int_{\{y:|y|\le z\}}K(z, 2\nu
s)dy)$$$$b_2(z;u(x-z,t-s), u(x-z,t-s))+(K(z,2\nu s){\bf
P}_z-$$$$\frac 1{4\pi|z|^3}(I-3e_ze_z^T) \int_{\{y:|y|\le z\}}K(y,
2\nu s)dy)g(x-z,t-s)\}dzds.\end{equation}

{\bf Theorem 2.4.} (\cite{Os}) {\em Let $h:R^3\to [0,\infty]$  and $\tilde h:R^3\to
[0,\infty]$ with $h$ locally integrable and $h,\tilde h$ jointly
satisfying
$$\int_{R^3}h^2(x-y)|y|^{-2} dy\le h(x)\mbox{and} \int_{R^3}\tilde h(x-y)|y|^{-1} dy\le h(x)
$$
for all $x\in R^3$ . If for all  $x\in R^3$ and $t>0$
$$(4\pi\nu t)^{-\frac 32}|\int_{R^3}u_0(x-y)
e^{-\frac{|y|^{2}}{4\nu t}} dy|\le \pi\nu \frac{h(x)}{11}$$ and
$$|g(x,t)|< (\pi\nu)^2 \frac{\tilde h(x)}{11}$$
then there exists a collection of probabilistic measures
$\{P_x:x\in R^3\}$ defined on a common measurable space $(\Omega,
{\cal F})$ and a measurable function $\Sigma:(0,\infty)\times
\Omega\to R^3$ such that
$$P_x(\{\omega:|\Sigma(t,\omega)|<\frac {2\pi \nu}{11}\mbox{ for all }
t>0)=1$$ for which  a weak solution $u(t,x)$ to the N-S can be
presented in the form
$$u(x,t)=h(x)\int_\Omega \Sigma(t,\omega)dP_x(\omega)\quad \mbox{ for all}
x\in R^3 , t>0.$$ Furthermore the solution $u$ is unique in the
class $$\{v\in (S'(R^3\times (0,\infty)))^3: |v(t,x)|<
\frac{2\pi\nu h(x)}{11}\} $$ for all $ \quad x\in R^3 , t>0.$}

Our survey is still far from being exhaustive. As already mentioned 
 the discussion of  the Euler-Lagrangian approach developed  by Constantin and Iyer  will be
 postponed to of the present paper the last section.

\section{A probabilistic representation of the solution to the Poisson equation}
\setcounter{equation}{0}

Within the framework of the approach developed in this paper we
intend to construct  diffusion processes associated with the
system (\ref{1.1}) (\ref{1.5}).  First we will start with
(\ref{1.5}) and recall some results concerning the solution of the
Poisson equation in an open domain  $G\subseteq R^3$.

First we recall that by the divergence theorem a $C^2(G)$ solution of
$-\Delta p=\gamma$ satisfies the integral identity
$$\int_G\nabla p\cdot\nabla \phi\,dx=-\int_G\gamma \phi \, dx$$
for all  $\phi\in C_0^1(G).$ In the space $W^{1,2}_0(G)$ which is the
completion of  $C_0^1(G)$ under the inner  product
$$\langle p,\phi\rangle=\int_G\nabla p\cdot\nabla\phi \,dx$$
the linear functional
$$F(\phi)=-\int_G \gamma\phi \,dx$$
may be extended to a bounded linear functional on the space
$W^{1,2}_0(G)$. Hence by the Riesz theorem there exists an element
$p\in W^{1,2}_0(G)$ satisfying $\langle p,\phi\rangle=F(\phi)$ for
all $\phi\in C_0^1(G).$ Then the existence of a
generalized solution to the Dirichlet problem $-\Delta p=\gamma$ and
$p=0$ on $\partial G$ is readily established. The question of
classical existence is accordingly transformed into the question of
regularity of generalized solution under the appropriately smooth
bounded conditions.

 We give in this section a brief summary
of a probabilistic approach to the solution of the Poisson
equation. We will try to give the probabilistic proofs of the necessary facts inasmuch  as they are  known. Proofs of similar statements can be found in \cite{BFR}.
The source for analytical results is the book by Gilbarg and Trudinger \cite{GTr}.

 Consider the
Poisson equation
\begin{equation}\label{3.1}
-\Delta p(x)=\gamma(x)
\end{equation} where $p$ and $\gamma$ are
scalar integrable functions defined on $ G$. A Newton potential
with density $\gamma$ is defined by
\begin{equation}\label{3.2}
N\gamma(x)=\frac{1}{4\pi}\int_{G}\frac{1}{\|x-y\|}\gamma(y)dy.\end{equation}
If $\gamma$ is regular and has a compact support then $N\gamma$ is
known to be a solution of the Poisson equation (\ref{2.1}).

To derive a probabilistic interpretation of the relation
(\ref{3.2}) we consider the generator ${\cal A}=\frac 12 \Delta $ of
a Wiener process $B(t)\in R^3$ defined on a given probability
space $(\Omega,{\cal F}, P)$. It is well known that on the space
$C_0(R^3)$ of all continuous functions vanishing at infinity the
Wiener process generates  the strongly continuous semigroup
$$T_t\gamma(x)=E\gamma(x+B(t)),\quad x\in R^3, t\ge 0, \gamma\in
C_0(R^3).$$

Given a function with a compact support in $G$ we extend it to the
whole space $R^3$ by zero.

By a direct computation we can check that
\begin{equation}\label{3.3}
\int_0^\infty E[\gamma
(x+B(t))]dt=\int_{R^3}\gamma(x+y)\int_0^\infty \frac{1}{(2\pi
t)^{\frac 32}}e^{-\frac{1}{2t}\|y\|^2}dtdy =$$
$$\int_{R^3}\frac{1}{2\pi\|y\|}\gamma(x+y)dy=2 N\gamma.\end{equation}

To prove that $p=2N\gamma$  solves the Poisson equation $$-\Delta
p=\gamma$$ we need some additional regularity properties of $
N\gamma$. 

{\bf Lemma 3.1.}{\em  Let $\gamma\in L^m(R^3)\cap L^q(R^3)$ with $1\le
m<3<q<\infty$. Then $N\gamma\in  C_0(R^3)$ and
$$\|N\gamma\|_\infty\le C_{m,q}(\|\gamma\|_m+\|\gamma\|_q).$$}

Proof.  First we note that for every  $l,m$ such that $\frac 1l+\frac
1r=1$  by H\"older inequality we have
\begin{equation}\label{3.4} E|\gamma(x+B(t))| =\frac {1}{(2\pi
t)^{\frac3 2}}
\int_{R^3}|\gamma(x+y)|e^{-\frac{\|y\|^2}{2t}}dy\le\end{equation}
$$
C_rt^{-\frac{ 3}{2}+\frac{ 3}{2l}}\|\gamma\|_r \le C_rt^{-\frac{
3}{2r}}\|\gamma\|_m,$$ since $-\frac{ 3}{2}+\frac{ 3}{2l}=-\frac{
3}{2r}.$

Finally we rewrite the left hand side of (\ref{3.1}) as
$$\int_0^\infty E[\gamma
(x+B(t))]dt=\int_0^1 E[\gamma (x+B(t))]dt+\int_1^\infty E[\gamma
(x+B(t))]dt$$ and applying the estimate (\ref{3.4})  for $r=q$ and
$r=m$ we derive
$$\int_0^\infty E[\gamma
(x+B(t))]dt\le C(\|\gamma\|_m+\|\gamma\|_q),$$ with $C=max (C_m,
C_q)$.

 By Sobolev embeddings  it is  known \cite{GTr} that if $\gamma\in L^1(R^3)$
then $N\gamma\in C(R^3)$. To check that $N\gamma\in C_0(R^3)$ we
note that for any $R>0$ we can rewrite the left hand side of
(\ref{3.3}) as
\begin{equation}\label{3.5}
\int_0^\infty E[\gamma (x+B(t))]dt=\int_0^\infty E[\gamma
(x+B(t))I_{\|B(t)\|>R}]dt+\end{equation}
$$
\int_0^\infty E[\gamma
(x+B(t))I_{\|B(t)\|\le R}]dt.
$$

 Let us prove that the first term on the right hand side of
 (\ref{3.5}) converges to 0 uniformly in $x$ as $R\to \infty$ and
 the second term converges to 0 as $\|x\|\to \infty$  for each $R$.
For the first term we apply the estimate (\ref{3.4}) to derive
$$\sup_{x\in R^3}E[|\gamma
(x+B(t))|I_{\{\|B(t)\|>R\}}]\le $$$$C(\|\gamma\|_p
+\|\gamma\|_q)(t^{-\frac 3{2m}}I_{\{[1,\infty)\}}(t)+t^{-\frac
3{2q}}I_{\{[0,1)\}}(t))$$ and
$$
\sup_{x\in R^3}E[|\gamma (x+B(t))|I_{\{\|B(t)\|>R\}}]\le
Ct^{-\frac
3{2}}\|\gamma\|_m(\int_{\{\|y\|>R\}}e^{-\frac{\|y\|^2}{2t}}dy)^{\frac
1q}\to 0$$ as $R\to \infty$. To obtain the estimate for the second
term we apply (\ref{3.4}) once again and obtain
$$E[|\gamma
(x+B(t))|I_{\|B(t)\|\le R}]\le$$$$
Ct^{-\frac{3}{2m}}(\int_{R^3}|\gamma(y)|^pI_{\{\|y-x\|\le
R\}}dy)^{\frac 1m} I_{\{[1,\infty)\}}(t)+
$$
$$Ct^{-\frac{3}{2q}}(\int_{R^3}|\gamma(y)|^qI_{\{\|y-x\|\le
R\}} dy)^{\frac 1q}I_{[0,1)}(t)),
$$
that yields  after the integration in time that the second term on
the right hand side of (\ref{3.5}) converges to 0, since $\gamma\in
L^m(G)\cap L^q(G)$ and is zero outside $G$.

To study derivatives of $N\gamma$ we apply the Bismut-Elworthy-Li
formula
$$\nabla_{x_i}E[\gamma(x+B(t))]=\frac{1}{t}E[\gamma(x+B(t))B_i(t)]$$
that holds for a regular $\gamma$.

{\bf Lemma 3.2.}{\em Let $\gamma\in L^m(R^3)\cap L^q(^3)$ for some
$1\le m<\frac 32<3<q <\infty$. Then $\nabla N\gamma\in C_0(R^3)$
and for each $x\in R^3$
\begin{equation}\label{3.6}
2\nabla_{x_i}N\gamma(x)=\int_0^\infty \frac 1 t
E[\gamma(x+B(t))B_i(t)]dt,\quad i=1,2,3.
\end{equation}
Moreover
\begin{equation}\label{3.7}
\|\nabla N\gamma\|_\infty\le C_{mq} (\|\gamma\|_m+\|\gamma\|_q).
\end{equation}}
 
 Proof.  By the H\"older inequality
\begin{equation}\label{3.8}
\frac 1 t E|\gamma(x+B(t))B_i(t)|dt =\frac {C}{t^{\frac
52}}\int_{R^3}|\gamma(x+y)y_i|e^{-\frac{\|y\|^2}{2t}}dy\le
\end{equation}
$$ \frac{C}{t^{\frac
52}}\|\gamma\|_mt^{\frac 12+\frac {3}{2q}}=
C_m\|\gamma\|_mt^{-\frac 12-\frac 3{2q}}.$$

Finally to give the sup estimate for the second derivative of
$N\gamma$  one has to apply the Schauder estimates and the
Bismut-Elworthy-Li formula.

Let us recall two more useful results (see \cite{GTr} theorem 4.5) concerning the Newton potential.

 {\bf Lemma 3.3.}{\em  Let $\gamma\in
L^q(R^3)\cap C^\alpha_b(R^3) $ with $1\le q\le \frac 32.$  Then
$N\gamma\in C^{2,\alpha}_b(R^3)\cap C_0(R^3)$,
$$\|N\gamma\|_{C^{2,\alpha}_b(R^3)}\le
C(\|\gamma\|_{L^q(R^3)}+\|\gamma\|_{C^\alpha_b(G)})$$ and  $p=
2N\gamma$ is the unique solution of the Poisson equation $$-\Delta
p=\gamma$$ in $C_0(R^3)\cap C^2(R^3)$.}

{\bf Theorem 3.4.}{\em Let $N\gamma\in C^2_0(R^3), \gamma\in
C^2_0(R^3)$ satisfy the Poisson equation $\Delta N\gamma=\gamma$ in
$R^3$. Then $N\gamma\in R^3$ and if $B=B_R(x_0)$ is any ball
containing the support of $N\gamma$ then
\begin{equation}\label{3.9}
\|\nabla^2N\gamma\|_{0,\alpha;B}\le
C_\alpha\|\gamma\|_{0,\alpha;B},\quad \|N\gamma\|'_{1,B}\le
CR^2\|\gamma\|_{0,B}.
\end{equation}}
 In the sequel we will need as well $L_q$ type
estimates for the Newtonian potential.

{\bf Lemma 3.5.}{\em The operator $N$ maps  $L^q(R^3)$ into $L^q(R^3)$
and there exists a positive constant  $C$ such that
\begin{equation}\label{3.10}
\|N\gamma\|_{L^q(R^3)}\le C\|\gamma\|_{L^q(R^3)}
\end{equation}}
Proof. By H\"older inequality we have
$$|N\gamma|(x)=|\int_{G}\gamma(y)(\Gamma(x-y))^{\frac 1q}(\Gamma(x-y))^{1-\frac
1q}dy|\le$$$$ \{\int_{G}|\gamma(y)|^q\Gamma(x-y)dy\}^{\frac
1q}\{\int_{G}\Gamma(x-y)dy\}^{1-\frac 1q}\le$$$$C
\{\int_{G}|\gamma(y)|^q\Gamma(x-y)dy\}^{\frac 1q}.$$

Next we obtain by Fubini's theorem
$$\int_{R^3}|N\gamma|^p(x)dx\le
\int_{R^3}C^p\{\int_{R^3}|\gamma(y)|^p\Gamma(x-y)dy\}dx=$$$$C^p\int_{G}\int_{R^3}|\gamma(y)|^p\Gamma(x-y)dydx=
C^p\int_{R^3}|\gamma(y)|^p(\int_{R^3}\Gamma(x-y)dx)dy\le$$$$
C_1\int_{R^3}|\gamma(y)|^pdy.$$

Note that all above results in this section are valid if we consider a bounded domain $G\subset R^3$ instead of $R^3$.
To get  further regularity properties of the Newton potential we
need more auxiliary  results.

Define the distribution $\nu_\gamma(\lambda)$ of the function
$\gamma: G\to R^1$  by
\begin{equation}\label{3.11}\nu_\gamma(\lambda)=|\{ x\in G:|\gamma(x)|>\lambda\}|\end{equation}
 where $|G|$ denotes
the Lebesgue volume of the domain $G$.

{\bf Lemma 3.6.}{\em Assume that $\gamma\in L^q(G)$ for some $q>0$.
Then
$$\nu_\gamma(\lambda)\le \lambda^{-q}\int_G|\gamma(x)|dx,$$
$$\int_G|\gamma(x)|^qdx=p\int_0^\infty\lambda^{q-1}\nu_\gamma(\lambda)d\lambda.$$}

Proof. It is easy to check that
$$\int_G|\gamma(x)|^pdx\ge \int_{\{\gamma>\lambda\}}|\gamma(x)|^pdx\ge
\lambda^p|\{x:\gamma(x)>\lambda\}|=\lambda^p\nu_\gamma(\lambda).$$
If $p=1$ we can apply the Fubini theorem to change the order of
integration
$$\int_G|\gamma(x)|dx=\int_G\int_0^{|\gamma(x)|}dtdx=\int_0^\infty\int_G
I_{\{x\in
G:\gamma(x)>\lambda\}}dxd\lambda=$$$$\int_0^\infty\nu_f(\lambda)d\lambda.$$

For arbitrary $q$ we have
$$\nu_{\gamma^q}(\lambda)=|\{x:\gamma(x)>
\lambda^{\frac {1}{q}}\}|=\nu_\gamma(\lambda^{\frac 1q})$$ and
hence
$$ p\int_0^\infty \lambda^{q-1}\nu_\gamma(\lambda)d\lambda=\int_0^\infty
\nu_{\gamma^q}(\lambda^q)d(\lambda^q)=\int_G|\gamma(x)|^qdx.$$

{\bf Lemma 3.7.}{\em Let $\gamma\in  L^q(G)$ for some $1<q<\infty$.
Then $ N\gamma \in W^{1,q}(G)$ and
\begin{equation}\label{3.12}
\|\nabla^2 N\gamma\|_{L^q(G)}\le C(q,G)\|\gamma\|_{L^q(G)}
\end{equation}
Moreover for $q=2$ the equality
\begin{equation}\label{3.13}
\int_{R^3}\|\nabla^2 N\gamma\|^2(x)dx=\int_G\gamma^2(x)dx
\end{equation}
 holds.}

Proof. The proof of this fact is  based on the Calderon-Zygmund
technique of cube decomposition and estimates of the function
$\nu_\gamma(\lambda)$ of the form (\ref{2.5}).

Let $ \tilde K$ be a cube in $R^3,$ $\gamma\ge 0$ integrable, and
finally fix $\kappa>0$ such that
$$\frac{1}{|\tilde K|}\int_{\tilde K}\gamma(x)dx\le \kappa.$$
Bisect $\tilde K$ into $2^3$ equal ( in volume) subcubes. Let $Q$ be a
set of those subcubes $K$  for which
$\frac{1}{|K|}\int_K\gamma(x)dx> \kappa$. For each of the
remaining subcubes ( which do not belong to $Q$) we repeat the
same procedure, that is bisect each one into $2^3$ sub-cubes  and
add those smaller ones, where $f$ is highly concentrated to $Q$.
Now repeating the procedure again and again we obtain a partition
of $\tilde K$. For any $K$ in $Q$ denote by $\hat K$ its
immediate predecessor.  Since $K\in Q,$ while $\hat K\notin Q$, we
have
$$\lambda<\frac{1}{|K|}\int_K\gamma(x)dx<\frac{1}{|K|}\int_{\hat K}\gamma(x)dx=
\frac{|\hat  K|}{|K|}\frac{1}{|\hat K|}\int_{\hat
K}\gamma(x)dx<2^3\lambda.$$

Set $F=\cup_{K\in Q}K, J=\tilde K\backslash F=\cap_{K\in
Q}K^C$. Note that each point in $J$  belongs to infinitely many
nested cubes with bounded concentration of $\gamma$ with diameters
converging to zero, that is $\frac{1}{|K_i|}\int_{\hat
K_i}\gamma(x)dx\le \kappa$, with $|K_i|\to 0.$ By the Lebesgue theorem
we deduce that $\frac{1}{|K_i|}\int_{\hat K_i}\gamma(x)dx\to
\gamma$ a.e. with respect to the  Lebesgue  measure, that is
$\gamma\le \kappa$ a.e. on $J$. Then we have  an
average estimate on $F$  and  a point-wise estimate on $J$.

At the second step we need the  Marcinkiewicz interpolation
theorem.

{\bf Marcinkiewicz interpolation theorem.}{\em Let $1\le q< r
<\infty$ and let ${\cal T}:L^q(G)\cap L^r(G)\to L^q(G)\cap L^r(G)$ be a
linear map. Suppose there exist constants $C_1, C_2$ such that
$\forall \gamma\in L^q(G)\cap L^r(G)$ and for any $\lambda>0$
$$
\nu_{{\cal T}\gamma}(\lambda)\le
\left(\frac{C_1\|\gamma\|_{L^q(G)}}{\lambda}\right)^q, \quad
\nu_{Q\gamma}(\lambda)\le
\left(\frac{C_2\|\gamma\|_{L^r(G)}}{\lambda}\right)^r.$$ Then for
any exponent $m$ such that  $q<m<r$ the map $\cal T$ can be extended to
a map from $ L^m(G) $ to $L^m(G)$  and
$$\|{\cal T}\gamma\|_{L^m(G)}\le KC_1^\alpha C_2^{1-\alpha}\|\gamma\|_{L^m(G)}.$$
 all $\gamma\in  L^q(G)\cap L^p(G)$ where $\frac 1m=\frac \alpha
 q+\frac {1-\alpha}r$ and the constant $K$ depends only on $m,q$ and
 $r$.}

At the end we define an operator ${\cal T}:L^2(G)\to L^2(G)$ by
${\cal T}\gamma=\nabla_i\nabla_jN\gamma$ to obtain the necessary result.

{\bf Theorem 3.8.(Calderon-Zygmund inequality)}{\em Let $\gamma\in
L^p(G)$, $1<q<\infty$ . Then the Newton potential $N\gamma=p\in
W^{2,q}(G)$, solves the Poisson equation  $\Delta p=\gamma$ a.e.
and
\begin{equation}\label{3.14}
\|\nabla^2 p\|_{L^q(G)}\le C \|\gamma\|_{L^q(G)},
\end{equation}
where $C$ depends only on $d$ and $q$. Furthermore , when $q=2$ we
have
$$\int_{R^3}\|\nabla^2 N\gamma(x)\|^2dx = \int_G\gamma^2(x)dx.$$}

For the proof of the above interpolation theorem and  theorem 3.8
see,  e.g., \cite {GTr}.

\section{ Probabilistic representations of \\ weak solutions of parabolic equations}
\setcounter{equation}{0}

In this section we adapt the results of the Kunita theory of
stochastic flows acting on Schwartz distributions
\cite{Ku1},\cite{Ku2} to the case under consideration.  The
considerations in this section are similar to  \cite{BW}.

Unlike the Kunita case we assume here that the coefficients of SDEs
under consideration are at most $C^{1+\alpha}$ -smooth with
$0<\alpha<1$, but on the other hand it is enough for our present
purpose to restrict ourself to  nonsingular initial data for the
Cauchy problem for parabolic equations and hence we consider
stochastic flows in Sobolev spaces ${\cal H}^k$ for $k=1,-1$.

It is more convenient for computational reasons to use sometimes
the Stratonovich form of the Ito equation. Recall that a process
$\xi(t)$ having  the Ito differential of the form
$$d\xi(t)=[a(\xi(t))+\frac 12 Tr\nabla \sigma(\xi(t))\sigma(\xi(t))]dt+\sigma(\xi(t))dw$$
has the Stratonovich differential of the form
$$d^S\xi(t)=a(\xi(t))dt+\sigma(\xi(t))\circ dw.$$

We say that condition {\bf C 4.1} holds if for all $t\in [0,T]$
the functions $g(t)$  and $\sigma$ belongs respectively to $
C^{1+\alpha}_b$ and  $C^{2+\alpha}_b$.

Throughout this section we assume that {\bf C4.1} holds. We shall
first give a brief review of the results which will be needed in the
sequel.

  Consider
a stochastic differential equation in the Stratonovich form
\begin{equation}\label{4.1}
d\xi(\tau)=-g(t-\tau,\xi(\tau))d\tau-\sigma(\xi(\tau))\circ
dw(\tau), \quad \xi(s)=x\in R^3,
\end{equation}
 $ 0\le s\le\tau\le t.$ Here $g(t,x)\in R^3$, $\sigma(t,x)\in
 R^3\times R^3$ and $w(t)\in R^3$ is a Wiener process.

Assuming that $g(t)\in {\bf C^1}(R^3)$ and  $\sigma(t)$ is a
$C^2$-smooth matrix  we are in the framework of the Kunita theory
\cite{Ku}  and
 know that there exists
a local $C^1$-diffeomorphism of $R^3$ generated by the  solution
$\xi_{s,x}(\tau)$ of (\ref{4.1}).

  Namely, by general results on the SDE theory the existence
and uniqueness of the solution $\xi_{s,x}(\tau)$ to (\ref{4.1})
are granted for a $C^1$- smooth bounded function $g$. Moreover, in
this case, one can prove that the solution $\xi^g_{s,x}(\tau)$ of
(\ref{4.1}) has a modification $\phi^g_{s,\tau}(x,\omega)$ such
that for all $\omega$ outside a null set ${\cal N}\subset \Omega $

1)$\phi^g_{s,\tau}(x,\omega)$ is continuous in   $(s,\tau,x)$, and
differentiable in $x$;

2)
$\phi^g_{\tau,\tau}(\phi^g_{s,\tau}(x,\omega),\omega)=\phi^g_{s,\tau}(x,\omega)$,
if $0<s<\tau< t$;

3) the mapping  $\phi^g_{s,\tau}(\omega):R^3\to R^3$  is a  $C^1$-
   diffeomorphism in  $R^3$.

The map    $\phi^g_{s,\tau}(\omega)$ is called a stochastic flow
of $C^1$- diffeomorphisms in  $R^3$.

 We will denote by  $(\phi^g_{s,\tau})^{-1}(\omega)=\psi^g_{\tau,s}(\omega)$
  the map inverse to the stochastic flow
 $\phi^g_{s,\tau}(\omega) $ and will write simply
$\psi^g_{\tau,s}(x)$ for $\psi^g_{\tau,s}(x,\omega)$. We check a
simple property  of an inverse stochastic flow.

{\bf Lemma 4.1.}{\em Consider the $\sigma$-algebras
$${\cal F}_s^w=\sigma\{w(\theta):\theta\in[0,s]\}\quad
 {\cal F}_{t,s}^{\hat w}=\sigma\{\hat w(\tau)-\hat w(\tau_1):s\le\tau_1\le \tau\le t\}$$
and a continuous bounded process $m(s)$ adapted to ${\cal F}_s^w$.
Then the process $f(s)=g(t-s)$ for $s\in [0,t]$ is   ${\cal
F}_{t,s}^{\hat w}$  adapted and for all $\alpha,\beta$ such that
$0\le\alpha\le\beta\le t$   we have
$$\int_\alpha^\beta f(\tau)dw(\tau)=\int_{t-\alpha}^{t-\beta}g(s)d\hat w(s).$$}

Proof. Note that since $ \hat w(s)=w(t-s)-w(t)$ we have $\hat
w(\beta)-\hat w(\alpha)=w(t-\beta)-w(t-\alpha),$  that yields
${\cal F}_{t-s}^w={\cal F}_{t,s}^{\hat w}$.

Now we consider a partition of the interval $[0,t]$
$$\{0= t_0\le
t_1\le \dots \le t_k\le t_{k+1}\le\dots\le t_N=t\}$$ such that
$|t_{k+1}-t_k|\to 0$ as $N\to\infty.$  Set  $\theta_k=t-t_k$ for
$k=1,\dots, N,$ then
$$\int_\alpha^\beta
f(s)dw(s)=\lim_{n\to\infty}\sum_{k=1}^Nf(t_k)[w(t_{k+1})-w(t_k)]=$$$$
\lim_{n\to\infty}\sum_{k=1}^Ng(t-\theta_k)[w(t-\theta_{k+1})-w(t-\theta_k)]=$$
$$-\lim_{n\to\infty}\sum_{k=1}^Ng(\theta_k)[\hat w(\theta_{k+1})-\hat w(\theta_{k})]=
-\int^{t-\alpha}_{t-\beta}g(s)d\hat w(s).$$

The main point of Kunita's theory is that the stochastic flow is a
bijection and that the inverse stochastic flow satisfies a couple of
SDEs which will be used for different purposes.  One of these SDEs
is given by the following lemma due to Malliavin ( see (\cite{IW},
lemma 5.2.2) or \cite{BD2}).

 {\bf Lemma 4.2.}{\em   Let $ \xi^g(\tau,x , w)$ be a
solution of the stochastic equation (\ref{4.1}) with $s=0$ or
equivalently of the SDE
  \begin{equation}\label{4.2} d \xi(\tau)=-g(t-\tau,
 \xi(\tau))d\tau+m(
 \xi(\tau))d\tau-\sigma( \xi(\tau)) d w(\tau) ,\quad \xi(s)=x
 \end{equation}
 where $m(x)=\frac 12 Tr \nabla\sigma(x)\sigma(x)$.
 Then, for every fixed $T>0$  we have
$$\xi(t-\theta,x, w)=\hat\xi(\theta,\xi(t,x,w),\hat w)$$
 for every $0\le \theta\le t$, and $x$, a.s. $(P^w)$.
}

In what follows  we need as well some generalizations of the
It$\hat{o}$ formula. The first one called the It$\hat{o}$ -Wentzel
formula reads as follows.

{\bf Lemma 4.3. (It$\hat{o}$-Wentzel formula)} {\em Assume that
the process $\xi(t)\in R^3$ has a stochastic differential of the
form
$$d\xi(t)=g(t,\xi(t))dt+\sigma(\xi(t))dw(t)$$
and the process $f(t,x)\in  R^3$ has a stochastic differential
$$df(t,x)=\Psi(t,x)dt+\Phi(t,x)dw(t)$$
with the same Wiener process $w(t)$. Let the vector field
$\Psi(t,x)\in R^3$ and the operator field $\Phi(t,x)\in R^3\times
R^3$  be $C^2$ smooth in $x$ and continuous in $t$. Then the
process $\eta(t)=f(t,\xi(t))$ has a stochastic differential
\begin{equation}\label{4.3}
df_m(t,\xi(t))=\Psi_m(t,\xi(t))dt+\Phi_{mk}(t,\xi(t))dw_k+ \nabla_i
f_m(t,\xi(t)) d\xi_i(t)+\end{equation}
$$
\frac 12  \nabla_i\nabla_jf_m(t,\xi(t))\sigma_{ik}(
\xi(t))\sigma_{jk}(\xi(t))dt
+\nabla_i\Phi_{mk}(t,\xi(t))\sigma_{ik}(\xi(t))dt.
$$}

{\bf Remark 4.4.} Note that (\ref{4.3}) can be rewritten in the
Stratonovich form as follows
\begin{equation}\label{4.5}
df_m(t,\xi(t))=\Psi_m(t,\xi(t))dt+\Phi_{mk}(t,\xi(t))dw_k+
\nabla_i f_m(t,\xi(t))\circ d\xi_i(t)+\end{equation}
$$
 +\nabla_i\Phi_{mk}(t,\xi(t))\sigma_{ik}(\xi(t))dt.
$$

We apply lemma 4.3 to check that  the inverse flow $\psi^g_{t,0}$
to the flow $\phi^g(0,t)$ (generated by the solution of the equation in
(\ref{4.1})) can also be represented as a solution of the
following stochastic equation
 \begin{equation}\label{4.4} d\psi^g_{t,0}(x)=\nabla
 \phi^g_{0,t}(\psi^g_{t,0})^{-1}g(t,x)dt+\nabla
 \phi^g_{0,t}(\psi^g_{t,0})^{-1}\sigma(x)\circ dw,
 \end{equation}
where $(\nabla\phi^g_{0,t})^{-1}$ is the inverse matrix of the
Jacobian matrix $\nabla\phi^g_{0,t}$ of the map
$\phi^g_{0,t}$.

Namely, we have the following statement  proved by Kunita (see
\cite{Ku} Theorem 4.2.2) in a slightly different context.

{\bf Theorem 4.5.}{\em Let {\bf C4.1} hold and  $\phi^g_{0,t}$ be
the solution of the equation (\ref{4.1}). Then the inverse flow
$[\phi^g_{0,t}]^{-1}=\psi^g_{t,0}$ satisfies  (\ref{4.4}).}

Proof. To verify the statement of the theorem note first that
the Jacobian matrix $\kappa^g(t)=\nabla \phi^g_{0,t}$ solves the
Cauchy problem for the stochastic equation
\begin{equation}\label{4.6} d\kappa(\tau)=-\nabla g(\tau,
\phi^g_{0,\tau}(y))\kappa(\tau)d\tau- \nabla
\sigma(\phi^g_{0,\tau}(y))\kappa(\tau)\circ d w(\tau) , \quad
\kappa(0)=I.
\end{equation}
Then, consider the stochastic process
 $$G(y,t)= \int_0^t\nabla
 \phi^g_{0,\tau}(y)^{-1}g(\tau,\phi^g_{0,\tau}(y))d\tau+$$$$\int_0^t\nabla
 \phi^g_{0,\tau}(y)^{-1}\sigma(\phi^g_{0,\tau}(y))\circ dw(\tau),$$
 and  compute $\phi^g_{0,t}(\psi^g_{t,0}(x))$,
 where $\psi^g_{t,0}$ has the  stochastic differential
 $$d\psi^g_{t,0}=dG(\psi^g_{t,0}(x),t).$$
 Set $\phi^g_{0,t}(y)= \phi^g(y,t).$
 By the It$\hat{o}$-Wentzell
 formula   we have
 $$\phi^g_{0,t}(\psi^g_{t,0}(x))= x+\int_0^td^S\phi^g(\psi^g_{\theta,0}(x),
 \theta)
 +\int_0^t\nabla\phi^g(\psi^g_{\theta,0}(x),\theta)\circ d\psi^g_{\theta,0}(x)=$$
 $$
x-\int_0^tg(\theta,\phi^g_{0,\theta}(\psi^g_{\theta,0}(x),\theta))d\theta-
\int_0^t\sigma(\phi^g_{0,\theta}(\psi^g_{\theta,0}(x),\theta))\circ
dw(\theta)+$$
$$
\int_0^t\nabla\phi^g(\psi^g_{\theta,0}(x),\theta)[\nabla\phi^g(\psi^g_{\theta,0}(
x),\theta)]^{-1}g(\theta,\phi^g_{0,\theta}(\psi^g_{\theta,0}(x)))d\theta+$$$$
\int_0^t\nabla\phi^g(\psi^g_{\theta,0}(x),\theta)[\nabla\phi^g(\psi^g_{\theta,0}(x),\theta)]^{-1}
\sigma(\phi^g_{0,\theta}(\psi^g_{\theta,0}(x),\theta))\circ
dw(\theta) =x.$$
 Hence,  $\phi^g_{0,t}(\psi^g_{t,0}(x))=x$ and thus $\psi^g_{t,0}$
 is the inverse to $\phi^g_{0,t}$.
\vskip 1mm

{\bf Remark 4.6.}  Recall that by lemma 4.2 the process
$\psi_{t,0}(x)=\hat\xi(t)$  along with (\ref{4.4}) satisfies the SDE
\begin{equation}\label{4.7}
d\hat\xi(\theta)=g(\theta,\hat\xi(\theta)))d\theta+\sigma(\hat\xi(\theta))\circ
d\hat w(\theta).
\end{equation}

Denote by $J^g_{0,t}(\omega)$ the Jacobian of the map
$\phi^g_{0,t}(\omega)$. Given $h\in{\cal H}^1$ and $f\in {\cal
H}^{-1}$ one can define the composition of $f$ with the stochastic
flow $\psi^g_{t,0}(x)$  as a random variable with values in ${\cal
H}^{-1}$ defined by the relation
\begin{equation}\label{4.8}
\langle S^g_{t,0}(\omega), h\rangle=\langle f,
h\circ\phi^g_{0,t}(\omega)J^g_{0,t}(\omega)\rangle,\quad h\in{\cal
H}^1,
\end{equation}
for any $t$ and $\omega\notin {\cal N}$. Note that if $\tilde f$
is a distribution of the form $\tilde f=f(x)dx$ where $f$ is a
continuous function then $\tilde f\circ\psi_{t,0}$ is just the
composition of the function $f$ with the map $\psi_{t,0}(\omega)$
and
$$\int_{R^3}f(\psi_{t,0}(y,\omega))h(y)dy=\int_{R^3}f(x)h(\phi_{0,t}(x,\omega))J_{0,t}(x,\omega)dx$$
by the formula of the change of variables.

{\bf Remark 4.7.} Consider the case of constant diffusion
coefficient $\sigma(x)\equiv \sigma$ and assume that the drift
possessed the property $div g=0$. Then (\ref{4.8}) has the form
\begin{equation}\label{4.9}
\langle S^g_{t,0}(\omega), h\rangle=\langle f,
h\circ\phi^g_{0,t}(\omega)\rangle,\quad h\in{\cal H}^1,
\end{equation}
since in this case $J^g_{0,t}(\omega)=Id$ is the identity map.

Consider  a linear PDE
\begin{equation}\label{4.10} \frac{ df}{dt}=L^gf-\gamma(t), \quad
f(0)=f_0,
\end{equation}
where
$$L^g=(g,\nabla)+L_0,$$
$$L_0f=\frac 12 F_{ij}\nabla_i \nabla_j f+ m_j \nabla_jf,$$
and
$$F_{ij}=\sigma_{ik}\sigma_{jk},\quad m_j=\nabla_j\sigma_{ik}\sigma_{jk},
 $$  and $g\in C^1$ is a given
bounded smooth function and  $f_0\in C^1$ (or more generally
$f_0\in{\cal D}'$).

 To construct a probabilistic representation of a weak
solution to (\ref{4.8}) in the case when the initial data  $f_0$
is a $C^1$ function (or even a distribution $f_0\in {\cal D}')$ we
consider the  stochastic process
\begin{equation}\label{4.11}
\lambda(t)=f_0-\int_0^t\gamma(\tau)\circ\phi^g_{0,\tau}d\tau
\end{equation}
and define its  composition  with a stochastic flow
$\psi^g_{t,0}(x)$  solving (\ref{4.7}). Recall that
$\psi^g_{t,0}(x)$  is inverse to the stochastic flow
$\phi^g_{0,t}(x)$ generated by the solution $\xi(t)$ of
(\ref{4.2}).

It is proved in \cite{Ku1} that the generalized solution of
(\ref{4.8}) is given by the generalized expectation of
$\lambda(t)\circ\psi^g_{t,0}(x)$.

To define the generalized expectation we consider the Sobolev
spaces $W^{k,q}$ or the weighted Sobolev spaces ${\cal S}_{k,q}$
defined in section 1 and check that
$E\lambda(t)\circ\psi^g_{t,0}(x)$ is well defined.

{\bf Lemma 4.8. }{\em  For each integer $k$ and $q>1$, $T>0$ there
are exist positive constants $c_{k,,q,T}, c'_{k,,q,T}$ depending
only on the flow $\psi_{t,0}$ such that  for any $t\in [0,T]$
\begin{equation}\label{4.12}
E\|\lambda(t)\circ\psi_{t,0}\|^q_k\le
c_{k,,q,T}\|f_0\|^q_{k,q}+c'_{k,,q,T}\int_0^t\|\gamma(\tau)\|^q_{k,q}d\tau,
\end{equation}
for all $f_0\in W^{k,q}$ and $\gamma(t)\in  W^{k,q}$.}

If $f_0\in  W^{k,q}$ then by this lemma for any $h\in W^{-k,q}$
there exists
$$\langle S_{0,t},h\rangle=E\langle
\lambda(t)\circ\psi_{t,0},h\rangle,$$ and $S_{0,t}$ can be considered
as an element from $W^{k,q}$. This element will be  called the
generalized expectation of $\lambda(t)\circ\psi_{t,0}$ and denoted
by $E[\lambda(t)\circ\psi_{t,0}]$.

For $k=1, q=2$ we consider $\langle \bar S^g_{t,0},h\rangle =
E\langle f\circ\psi^g_{t,0}, h\rangle$ which is a continuous
linear functional on ${\cal H}^1$ and can be regarded as an
element of ${\cal H}^{-1}$ .  Set
\begin{equation}\label{4.13}
U_{t,0}f=E[f\circ\psi_{t,0}]
\end{equation}
and call it the generalized expectation of
$f\circ\psi_{t,0}$.
It is easy to see that  $U_{t,0}$ is a linear map form ${\cal H}^{-1}$ into itself.
Moreover the family $U_{t,s}f=E[f\circ\psi_{t,s}]$ possesses the
evolution property $U_{t,\tau}U_{\tau,s}=U_{t,s}$ for any $0\le
s\le t\le T$.  It can be immediately deduced from the evolution
properties of $\phi^g_{s,t}$ and $J^g_{s,t}$.

Finally we  compute the infinitesimal operator of the evolution
family $U_{t,s}$. To this end we need  a version of the
It$\hat{o}$ formula.

{\bf Theorem 4.9.}(The generalized It$\hat{o}$ formula) {\em
 Let
$f(t)\in {\cal H}^1$
 be a continuous  in $t$ nonrandom function. Then, given
 stochastic flows $\phi_{0,t}, \psi_{t,0}$ generated by (\ref{4.1}), (\ref{4.6})
the following relations hold
$$f(t)\circ\phi^g_{0,t}=f(0)+
\int_{0}^{t}[\frac{\partial f(\theta)}{\partial\theta} \circ
\phi^g_{0,\theta}+L_0f (\theta) \circ \phi^g_{0,\theta}]d\theta+
$$
$$
  \int_{0}^{t} \nabla_if
(\theta)\circ\phi^g_{0,\theta} d[\phi^g]^i_{0,\theta}$$  and
\begin{equation}\label{4.14}
f(t)\circ\psi^g_{t,0}=f(0)+\int_{0}^{t}[\frac{\partial }{\partial
\theta}[f(\theta)]\circ\psi^g_{\theta,0} + L_0[f(\theta)\circ
\psi^g_{\theta,0}]]d\theta+\end{equation}$$\int_{0}^{t} \nabla_i[
f(\theta)\circ\psi^g_{\theta,0}] \sigma dw+ \int_{0}^{t} \nabla_i[
f(\theta)\circ\psi^g_{\theta,0}] g(\theta) d\theta.
$$
Here   we understand the action of the operator $L_0$
 in the sense of
generalized functions.}

 The proof of  theorem 4.9 employs the classical It$\hat{o}$ formula for $C^2$- smooth functions $f_\varepsilon$
 that approximate the $C^1$  function
 $f$, uses  equations (\ref{4.2}) and (\ref{4.6}) for the flows
 $\phi_{0,t}$ and $\psi_{t,0}$,  respectively, and  then justifies
 the passage to the limit under the integral sign in the integral identity.
  The details can be found in \cite{Ku1} for a much more general case.

Let us come back to the parabolic equation (\ref{4.10}) and set
$\gamma=0$.  We can show that the stochastic flow $\psi_{t,0}$
gives rise to an evolution family acting in spaces of of
distributions and the function $f(t)=E[f_0\circ\psi_{t,0}]$ is a
weak solution of (\ref{4.10}) with $\gamma=0$.

{\bf Theorem 4.10. }{\em  Assume that the coefficients of the
stochastic equation  (\ref{4.2}) satisfy {\bf C 4.1} and
$\psi^g_{t,0}$ is generated by the solution of (\ref{4.6}). Then,
for any functions
 $f_0,g(t)\in {\cal H}^1$  the relation
 $$f(t)=E[f_0\circ \psi^g_{t,0}]$$
defines the unique generalized solution to the problem (\ref{4.10})
with $\gamma=0$. The restriction of  $U^g_{t,0}$ to ${\cal H}^1$
defines a strongly continuous family of evolution mappings acting on
the space ${\cal H}^1$. The domain of definition of its
infinitesimal operator
  ${\cal A}^g$  (in a weak sense) contains the subspace
 ${\cal H}^1$ and ${\cal A}^gf=L^gf$ for any  $u\in
{\cal H}^1$.} \vskip 1mm

 Proof.
   From the relation (\ref{4.13}) and the properties of stochastic flows we  deduce
that the relation
$$  \langle U^g(t)f_0,h \rangle
  =  \langle E[f_0\circ\psi^g_{t,0}], h \rangle  $$
 defines a continuous linear functional
on ${\cal H}^1$. Thus, we can treat $U^g(t)f_0$ as an element from
${\cal H}^1$. It follows from the representation
$U^g_{t,0}f_0(x)=E[f_0\circ \psi^g_{t,0}(x)] $ that $U^g(t)$ is a
linear mapping from the space ${\cal H}^1$ into itself. Note that
the above definition of the family $U^g_{t,s} $ through the
integral identity allows to check that it possesses the evolution
property
$$\langle
U^g_{t,\tau}U^g_{\tau,s}f_0, h\rangle=\langle U^g_{t,s}f_0,
h\rangle.$$

  Indeed, by the Markov property of the process
  $\psi^g_{t,0}(x)$ we deduce that
$$\langle U^g_{t,s}f_0, h\rangle=\langle U^g_{t,\tau}U^g_{\tau,s}f_0,\, h\rangle=
E[\langle U^g_{\tau,s}f_0,\, h\circ \phi^g_{\tau,t}
J^g_{\tau,t}\rangle]=
$$
$$ E[\langle f_0,\,[v\circ \phi^g_{s,\tau}]J^g_{s,\tau}
\rangle|_{v=h\circ\phi^g_{\tau,t}J^g_{\tau,t}}]= E\langle f_0,\,
[h\circ\phi^g_{s,\tau}\circ\phi^g_{\tau,t}]
 J^g_{s,\tau}]\circ\phi^g_{\tau,s} J^g_{\tau,t}\rangle=$$
$$
E[\langle f_0,\, [h\circ\phi^g_{s,t}]J^g_{s,t} \rangle]=\langle
U^g_{s,t}f, h\rangle.
$$
Now we  apply  the generalized It$\hat{o}$ formula to obtain the
relation
$$E[f_0\circ\psi^g_{t,0}]=f_0+E[\int_{0}^{t}L^g (f_0\circ
\psi^g_{\theta,0})d\theta].$$ Note that in the latter expression
each summand  belongs to ${\cal H}^1$.  In addition,
 $$E [\int_{0}^{t}\langle L^g(f_0\circ
\psi^g_{\theta,0}),\, h\rangle d\theta]=\int_{0}^{t}\langle
E[f_0\circ\psi^g_{\theta,0}],(L^g)^*h \rangle
d\theta=$$$$\int_{0}^{t}\langle L^g(E
[f_0\circ\psi^g_{\theta,0}]),h \rangle d\theta,$$ that yields
 $$E[f_0\circ \psi_{t,0}^g]=f_0+\int_{0}^{t}L^g(E[f_0\circ
\psi_{\theta,0}^g])d\theta.$$ In other words
$$U^g_{t,0}f_0=f_0+\int_0^tL^gU^g_{\theta,0}f_0d\theta.$$
As the result we get that
$f(t)=E[f_0\circ\psi_{t,0}]$ satisfies (\ref{4.8}) and
$f(0)=f_0$.

One can prove the corresponding result in the case $\gamma(t)\neq 0$
in a similar way  applying the above reasons to $\lambda(t)$ of the
form (\ref{4.11}) instead of $f_0$.

{\bf Theorem 4.11.} {\em Given tempered distributions  $f_0$ and
$\gamma(t)$ define $\lambda(t)$ by (\ref{4.11}). Then
$U(t)=E[\lambda(t)\circ\psi^g_{t,0}]$ defines  the unique solution
of equation (\ref{4.10}) if $\psi^g_{t,0}$ satisfies (\ref{4.6}) and
$\phi_{0,t}$ is its inverse.}

Proof.  By the generalized Ito formula we get
$$ \lambda(t)\circ\psi^g_{t,0}=f_0-\int_0^t
\gamma(\tau)d\tau+\int_0^t\nabla_i(\lambda(\tau)\circ\psi^g_{\tau,0})g(\tau)d\tau
+$$
$$\int_0^t\nabla_i(\lambda(\tau)\circ\psi^g_{\tau,0})\sigma(\tau)dw(\tau)+
\int_0^tL^g_0(\lambda(\tau)\circ\psi^g_{\tau,0})d\tau.$$ As a
consequence we get
$$\lambda(t)\circ\psi^g_{t,0}=f_0+\int_0^t\nabla_i(\lambda(\tau)
\circ\psi^g_{\tau,0})\sigma dw(\tau)+$$
\begin{equation}\label{4.15}\int_0^tL^g(\lambda(\tau)
\circ\psi^g_{\tau,0})d\tau-\int_0^t\gamma(\tau)d\tau.\end{equation}
Each term in (\ref{4.15}) has a generalized expectation as an
element of ${\cal S}'$. The generalized expectation of the second
term in the right hand side of (\ref{4.15})  is equal to zero. For
the third term we have
$$\langle E\left[\int_0^tL^g(\lambda(\tau)\circ\psi^g_{\tau,0})d\tau\right], h\rangle=
E\left[\int_0^t\langle\lambda(\tau)\circ\psi^g_{\tau,0},
[L^g]^*h\rangle d\tau\right] =$$
\begin{equation}\label{4.16}
\int_0^t \langle E[\lambda(\tau)\circ\psi^g_{\tau,0}],
[L^g]^*h\rangle d\tau=\int_0^t \langle
L^g(E[\lambda(\tau)\circ\psi^g_{\tau,0}]),h\rangle d\tau.
\end{equation}
Hence
\begin{equation}\label{4.17}
E[\lambda(t)\circ\psi^g_{t,0}]=f_0+\int_0^tL^g(E[\lambda(\tau)\circ\psi^g_{\tau,0}])d\tau-
\int_0^t\gamma(\tau)d\tau.
\end{equation}
Differentiating each term with respect to $t$ we check that
$U(t)=E[\lambda(t)\circ\psi^g_{t,0}]$ satisfies (\ref{4.10}). In
addition $\lim_{t\to 0} \langle U^g(t),h\rangle=\langle
f_0,h\rangle$, that is $\lim_{t\to 0} \langle U^g(t)=f_0$ and we
proved that $U^g(t)$ solves the Cauchy problem (\ref{4.10}).

 To prove the
uniqueness of the solution to (\ref{4.8}) suppose to the contrary
that there exist two solutions $f(t)$ and $\tilde f(t)$ to
(\ref{4.10}).  Then the function $u(t)=f(t)-\tilde f(t)$ satisfies
$\frac{d u(t)}{dt}=L^g u(t)$   and $\lim_{t\to 0}u(t)=0$. Fix $t$
and choose a function $h(t,\cdot)\in {\cal D}$. Then there exists a
solution $h(\tau,x), \, 0\le s\le \tau\le t, x\in R^3,$ to the
Cauchy problem
     $$\frac{\partial h(\tau,x)}{\partial
\tau}+[L^g]^*h(\tau,x)=0, \quad \lim_{\tau\to t}h(\tau,x)=h(t,x).$$
    If the coefficients   $a^g$ and  $ \sigma^g$ are $C^1$-smooth , then
     there exists a unique classical solution to this Cauchy problem. As a result,
    $$\langle  u(t), h\rangle =\int_{0}^{t}  \langle \frac{d}{d\theta}u(\theta),
    h(\theta)\rangle
    d\theta +\int_{0}^{t}\langle u(\theta),\frac{d}{d\theta}h(\theta \rangle d\theta=$$
$$\int_{0}^{t} \langle L^gu(\theta), h(\theta)\rangle  d\theta -\int_{0}^{t}
    \langle u(\theta), [L^g]^*h(\theta)\rangle d\theta=0.$$

\section{A probabilistic approach to the \\Navier-Stokes system}
\setcounter{equation}{0}

Let us come back to  the Navier-Stokes system

\begin{equation}\label{5.1}
 \frac{\partial u}{\partial t}+(u,\nabla )u=\frac{\sigma
^2}2\Delta u-\nabla p ,\quad u(0,x)=u_0(x),\quad x\in R^3
\end{equation}

\begin{equation}\label{5.2}
 -\Delta p=\gamma,\end{equation}
with $\gamma$ defined by (\ref{1.3}).

Our main purpose in this section is to construct a diffusion process
that allows us to obtain a probabilistic representation of a weak
solution to (\ref{5.1}), (\ref{5.2}). To be more precise we intend
to reduce the solution of this system to  solution of a certain
stochastic problem, to solve it and then to verify that in this way
we have constructed a weak solution of (\ref{5.1}), (\ref{5.2}).

 Let as above
$w(t), B(t)$ be  standard $R^3$-valued independent Wiener processes
defined on a probability space  $(\Omega,{\cal F},P)$.
 Given a bounded
measurable function $f(x)$ and a stochastic process $\xi(t)$ we
denote $E_{s,x}f(\xi(t))\equiv Ef(\xi_{s,x}(t))$ a conditional
expectation under the condition $\xi(s)=x$.

In section 2 we recalled the probabilistic approach  developed in
our previous paper \cite{AB} that allows to construct a
probabilistic representation of a $C^2$-smooth (classical)
solution to (\ref{5.1})-(\ref{5.2}) via the solution of the
stochastic problem
 \begin{equation}\label{5.4}d\xi(\tau)=-u(t-\tau,\xi(\tau))d\tau+
\sigma dw(\tau),\end{equation}
\begin{equation}\label{5.5}
u(t,x)=E_{0,x}[u_0(\xi(t))+\int_0^t\nabla p(t-\tau,\xi(\tau))d\tau]
\end{equation}
\begin{equation}\label{5.6}
-2p(t,x)=E[\int_0^\infty\gamma(t, x+B(t))dt]= E[\int_0^\infty{\rm
tr}[\nabla u]^2(t, x+B(t))dt].
\end{equation}

In this section  we consider a similar stochastic  system but
now we choose to invert the time direction of the stochastic
process itself rather then of the function $u$ to obtain the
possibility to reduce a construction of a generalized solution to
the Navier-Stokes system to the construction of a  solution of a
stochastic problem.

Our considerations will be based on the result of sections 3 and
4. Note that since we consider the case where the diffusion
coefficient $\sigma$ is constant the Ito form and the Stratonovich
form of a stochastic equation coincide.

 Let as above
$w(t), B(t)$ be  standard $R^3$-valued independent Wiener
processes defined on a probability space  $(\Omega,{\cal F},P)$.

Let $\phi_{0,t}(y)$ be a stochastic process satisfying the
stochastic equation
$$d\phi_{0,t}(y)=u(t,\phi_{0,t}(y))dt-\sigma dw(t),\quad \phi_{0,0}=y$$
and the stochastic process $\lambda(t)$ be of the form
\begin{equation}\label{5.07} \lambda(t)=u_0-\int_0^t\nabla p(\tau, \phi_{0,\tau})d\tau.\end{equation}
Consider the system
\begin{equation}\label{5.7} d\psi_{t,\theta,}(x)=-u(\theta,
 \psi_{t,\theta}(x)) d\theta + \sigma d\hat w(\theta), \quad \psi_{t,t}(x)=x,
 \end{equation}
\begin{equation}\label{5.8}
u(t,x)=E[u_0(\psi_{t,0}(x))-\int_0^t\nabla p(\tau,
\psi_{t,\tau}(x))d\tau].
\end{equation}

\begin{equation}\label{5.9}
-2\nabla p(t,x)=E[\int_0^\infty\frac 1\tau\gamma(t,
x+B(\tau))B(\tau)d\tau],
\end{equation}
 where $\gamma$ is given by (\ref{1.3}) and prove the existence and uniqueness of its solution.

To this end we apply the Picard principle  to the solution of the
stochastic system and construct a solution to
(\ref{5.7})-(\ref{5.9}) by the successive approximation technique.

Set
\begin{equation}\label{5.10}
u^1(t,x)=u_0(x), \quad \psi^0_{t,0}(x)=x, p^1(t,x)=0
 \end{equation}
and consider a family of stochastic processes
$\psi^k_{t,\theta}(x)$ and families  of vector fields $u^k(t,x)$
and scalar functions $p^k(t,x)$ given by the following relations

\begin{equation}\label{5.11}
d\psi^k_{t,\theta}=-u^k(\theta,\psi^k_{t,\theta})d\theta+\sigma
d\hat w(\theta),\quad  \psi^k_{t,t}=x,
 \end{equation}

\begin{equation}\label{5.12}
u^{k+1}(t,x)= E[u_0(\psi^k_{t,0}(x))-\int_0^t\nabla
p^{k+1}(\tau,\psi^k_{t,\tau}(x))d\tau],
 \end{equation}
\begin{equation}\label{5.13}
-2 p^{k+1}(t,x)=\int_0^\infty E[ \gamma^{k+1}(t,x+B(\tau))]d\tau,
 \end{equation}
where
\begin{equation}\label{5.14}
\gamma^{k+1}(t,x)=Tr[\nabla u^{k}(t,x) \nabla u^{k+1}(t,x)].
 \end{equation}

Note that  for a fixed $k$ the first stochastic equation
(\ref{5.11}) that determines the family of stochastic processes
$\psi^k_{t,0}(x)$ may be solved independently on equations
(\ref{5.12})-(\ref{5.14}).  Then given the process
$\psi^k_{t,0}(x)$ and
  keeping in mind the properties of the function $p^k$
 that satisfies the Poisson equation
\begin{equation}\label{5.15}
-\Delta p^{k+1}(t,x)=\gamma^{k+1}(t,x),
 \end{equation}
one has to compute $\nabla p^k(t,x)$, $u^{k+1}(t,x)$
 by  (\ref{5.12}), (\ref{5.13}).

 To investigate the convergence of the stochastic processes
$\psi^k_{t,0}(x) $ and functions $u^k(t,x),p^k(t,x)$ defined above
we need some auxiliary results concerning the behavior of
solutions of stochastic equations.

Let $g\in V$ be a given function. Consider the stochastic equation

\begin{equation}\label{5.16}
d\psi^g_{t,\theta}=-g(\theta)\circ\psi^g_{t,\theta}d\theta+\sigma
d\hat w(\theta), \quad \psi^g_{t,t}(x)=x
 \end{equation}
and define vector fields $u^g(t,x)$ and $\nabla p^g(t,x)$ by

\begin{equation}\label{5.17}
u^g(t,x)= E[u_0(\psi^g_{t,0}(x))-\int_0^t\nabla
p^g(\tau,\psi^g_{t,\tau}(x))d\tau],
 \end{equation}

\begin{equation}\label{5.18}
-2p^g(t,x)=\int_0^\infty E[ \gamma^g(t,x+B(\tau))]d\tau,
 \end{equation}

\begin{equation}\label{5.19}
\gamma^g(t,x)=Tr[\nabla g\nabla u^g](t,x).
 \end{equation}

 Recall that $p^g$ solves the Poisson equation
\begin{equation}\label{5.20}
-\Delta p^g=Tr[\nabla g\nabla u^g].
 \end{equation}
 To investigate the convergence of the stochastic processes
$\psi^k_{t,0}(x) $ and functions $u^k(t,x),p^k(t,x)$ defined above
we need some auxiliary results concerning the behavior of
solutions of stochastic equations. Moreover along with the system
(\ref{5.16}) -- (\ref{5.18}) we will need the system to describe
the process $\eta^k(\tau)=\nabla\psi^k_{t,\tau}(x)$ and the
functions $\nabla u^k(t,x)$ and $\nabla p^k(t,x)$.

To derive the necessary apriori estimates we start with the
consideration of a linearized system.

Let $g(t)$ be a given vector field. Consider the stochastic
equation

\begin{equation}\label{5.21}
d\psi^g_{t,\theta}=-g(\theta)\circ\psi^g_{t,\theta}d\theta+\sigma
d\hat w(\theta), \quad \psi^g_{t,t}(x)=x
 \end{equation}
and define the vector fields $u^g(t,x)$ and $\nabla p^g(t,x)$ by

\begin{equation}\label{5.22}
u^g(t,x)= E[u_0(\psi^g_{t,0}(x))-\int_0^t\nabla
p^g(\tau,\psi^g_{t,\tau}(x))d\tau],
  \end{equation}

\begin{equation}\label{5.23}
-2 p^g(t,x)=\int_0^\infty E[ \gamma^g(t,x+B(\tau))]d\tau,
 \end{equation}

\begin{equation}\label{5.24}
\gamma^g(t,x)=Tr[\nabla g\nabla u^g](t,x).
\end{equation}

 Finally we  derive  from (\ref{5.23})  the relation
\begin{equation}\label{5.25}
-2 \nabla p^g(t,x)=\int_0^\infty E[\frac 1\tau
\gamma^g(t,x+B(\tau))B(\tau)]d\tau
\end{equation} by applying the
Bismut -- Elworthy -- Li formula first  checking  the
conditions that validate such an application are satisfied.  Below we will need
some estimates of a solution of the Poisson equation  from section
3. For convenience of references we formulate them in the
following statement.

{\bf Lemma 5.1.} 

{\em 1. Let $\gamma^g\in L^q(R^3)\cap L^m(R^3)$
for some $1\le q<\frac 32<3<m<\infty$. Then
$$\|\nabla p^g\|_\infty\le C_{qm}(\|\gamma^g\|_q+\|\gamma^g\|_m)$$
$$\|\nabla_i\nabla_jp^g\|_\infty\le
C(\|\gamma^g\|_q+[\gamma^g]_\alpha).$$

2. Let  $\gamma^g\in L^r(R^3)$ for $1<r<\infty$. Then $p^g\in
W_0^{2,r}(R^3)$
 and the Calderon- Zygmund
inequality
$$\|\nabla_i\nabla_jp^g\|_{r,loc}\le C_1\|\gamma^g\|_{r,loc}$$
holds.}

\vskip 1mm

Let  $ {Lip}$ be the subspace of the space $C(R^1\times R^3,R^3)$
of continuous ( in $t\in [0,T], x\in R^3$), bounded functions
which
 consists of Lipschitz-continuous  (in $x$) functions
$g$ such that
$$\|g(t,x)-g(t,y)\|\le L_g(t)\|x-y\|, \quad t\in [0,T]\quad
x,y\in R^3,$$
 where $\|\cdot\|$ is the norm in $R^3$.

 Condition   {\bf C 5.1}

Let $g(t,x)\in R^3$  be  a vector field defined on $[0,T]\times
R^3$ that belongs to $C^{1,\alpha}(R^3,R^3), 0<\alpha\le 1 $ for a
fixed $t\in [0,T]$ and satisfies the following estimates:

 $ 1.\quad   \|g(t)\|_{ L^q_{loc}}\le N_g(t)$ for some $q$ to be specified below,
 $\|g(t)\|_\infty\le K_g(t)$  and
$$\|g(t,x)-g(t,y)\|\le L_g(t)\|x-y\|,\quad \|\nabla g(t,x)-\nabla g(t,y)\|\le L_g^1(t)\|x-y\|.$$

$2.\quad  \|\nabla g(t)\|_\infty\le K^1_g(t), \quad \|\nabla
g(t)\|_{r, loc}\le N^1_g(t),$ where $K_g(t), L_g(t)$, $N_g(t)$ and
$K^1_g(t), L^1_g(t), N^1_g(t)$ are positive functions bounded on
an interval $ [0,T]$, $r=m$ and $r=q$ for $1<q<\frac
32<3<m<\infty$.

Set $\psi(\tau)=\psi_{t,\tau}(x)$ and consider the stochastic
equation

 \begin{equation}\label{5.26}
  \psi(\tau)=x-\int_\tau^t
 g(\tau_1,
 \psi(\tau_1))d\tau_1 +\int_\tau^t\sigma d\hat w(\tau_1),
\end{equation}
with $0\le \tau\le t<T$ for a certain constant $T$.  If we are
interested in the particular dependence of the process
$\psi(\tau)$ on the parameters $t,x$ and $g$, we  write
$\psi(\tau)=\psi_{t,x}^g(\tau)$.

 {\bf Lemma 5.2.} {\em Assume that   {\bf C 5.1} holds. Then there exists
a unique solution  $\psi^g_{x}(\tau)$ of (\ref{5.21}) that
satisfies the following estimates:
\begin{equation}\label{5.27}
E\|\psi^g_x(\tau)\|^2\le 3[\|x\|^2 +\sigma^2(t-\tau)+
(t-\tau)\int_{\tau}^{t}[K^2_g(\tau_1)]d\tau_1] ,\end{equation}
\begin{equation}\label{5.28}
E\|\psi^g_{x}(\tau)- \psi^g_{y}(\tau)\|\le \|x-y\|
e^{\int_{\tau}^{t}L_g(\theta) d\theta}, \end{equation}
\begin{equation}\label{5.29}E\|\psi^g_{x}(\tau)- \psi^{g_1}_{x}(\tau)\|\le
 \int_{\tau}^{t}\|g(\tau_1)-g_1(\tau_1)\|_\infty d\tau_1
e^{\int_{\tau}^{t} L_g(\theta)d\theta}. \end{equation}}

 \vskip 1mm

Proof. The proof of the estimates of this lemma is standard.
 We only show  the proof of (\ref{5.28}).  In view of  {\bf
C 5.1} we have
$$E\|\psi_x^g(\tau)-\psi^g_y(\tau)\| \le
\|x-y\|+
\int_{\tau}^{t}L_g(\tau_1)\|\psi^g_x(\tau_1)-\psi^g_y(\tau_1)\|d\tau_1,$$
where $0\le \tau\le t\le T$ with some constant $T$ to be chosen
later.
 Finally, by Gronwall's lemma, we get
$$E\|\psi^g_x(\tau)-\psi^g_y(\tau)\| \le
\|x-y\|e^{\int_{\tau}^{t} L_g(\theta)d\theta}.$$

 Along with the equations
(\ref{5.21})-(\ref{5.23}) we will need below the equations for the
mean square derivative $\eta(t)=\nabla \psi_{t.0}(x)$ of the
diffusion process $\psi_{t,0}(x)$ that satisfies (\ref{5.21}) and
the gradient $v(t,x)=\nabla u(t,x)$ of the function $u(t,x)$ of
the form (\ref{5.22}).

{\bf Lemma 5.3} {\em Assume that  {\bf C 5.1} holds. Then the
process $\eta^g(\tau)=\nabla\psi^g_{t,\tau}$ satisfies the
stochastic equation
\begin{equation}\label{5.30}
d\eta^g(\tau)=-\nabla
g(\tau,\psi^g_{t,\tau})\eta^g(\tau)d\tau,\quad \eta^g(t)=I,
\end{equation} where $I$ is the identity map. Furthermore the
process $\eta^g(\tau)$ possesses the following properties.

The determinant  $det\,\eta(\tau)$ is equal to 1, i. e.
$$det \,\eta^g(\tau)=J_{t,\tau}=1,$$
and the estimate
 \begin{equation}\label{5.31}\|\eta^g(\tau)\|\le  e^{\int_\tau^t
K^1_g(\theta)d\theta}\end{equation} holds.

In addition the following integration by parts formula is valid
\begin{equation}\label{5.32}
\int_{R^3}f(\psi^g_{x}(\tau))dx=\int_{R^3}f(x)dx,\quad f\in
L^1(R^3). \end{equation}}

Proof. Under the above assumptions the first statement immediately
follows from the results of the stochastic differential equation
theory. By a direct computation one can check that $J_{t,\tau}$
satisfies the linear equation
$$dJ_{t,\tau}=-{ div}\, g(\psi^g_{t,\tau})J_{t,\tau}d\tau,\quad
J_{t,t}=I$$ and since $div\, g=0$ we get the second statement that
yields the integration by parts formula (\ref{5.32}). Finally
(\ref{5.31}) is deduced from the inequality
$$E\|\eta(\tau)\|\le 1+\int^t_\tau K^1_g(\theta)E\|\eta(\theta)\|d\theta$$
by the Gronwall lemma.

In the sequel we denote by $\eta^{x,g}(t)$ the solution of the
equation
$$d\eta^{x,g}(\tau)=\nabla
g(t,\psi_{t,\tau}(x))\eta^{x,g}(\tau)d\tau,\quad \eta^{x,g}(0)=I$$
if we will be interested in the properties of the process
$\eta^{x,g}(t)$.  One can easily check that
$$\|\eta^{x,g}(\tau)-\eta^{y,g}(\tau)\|\le \int_\tau^t\|\nabla
g(\theta,\psi_{t,\theta}(x))- \nabla
g(\theta,\psi_{t,\theta}(y))\|d\theta
e^{\int_\tau^tK^1_g(\theta)d\theta}$$$$ \le \int_\tau^t
L^1_g(\theta)\|\psi_{t,\theta}(x))-\psi_{t,\theta}(y))\|d\theta$$
and by (\ref{5.28}) we have
$$\|\eta^{x,g}(\tau)-\eta^{y,g}(\tau)\|\le C(\tau)\|x-y\| $$
where $C(\tau)$ is a bounded function over a certain interval $[0,
T_1]$ depending on $g$.

Let us state conditions on initial data $u_0$ of the N-S system.

We say that {\bf C 5.2} holds if  for $0<\alpha\le 1$ the initial
vector field $u_0 \in C^{1,\alpha}$ satisfies the following
estimates
$$\|u_0\|_\infty\le K_0, \quad \|\nabla u_0\|_\infty\le K^1_0,\quad
\|u_0\|_{r,loc}\le M_0,\quad  \|\nabla u_0\|_{r,loc}\le M_0^1$$
 with $r$ to be specified below and let $L_0, L^1_0$ be Lipschitz constants for the functions $u_0$ and $\nabla u_0$ respectively.

{\bf Lemma 5.4.} Assume that  $g(t,x)$ satisfies {\bf C 5.1} and
$u_0$ satisfies {\bf C 5.2} with $r=q$ and $r=m$ for
$1<q<\frac32<3<m<\infty$. Then  the vector field $u^g(t,x)$ given
by
\begin{equation}\label{5.33}
u^g(t,x)=E[u_0(\psi^g_{t,x}(0))-\int_0^t\nabla
p^g(\tau,\psi^g_{t,\tau}(x))d\tau] \end{equation} satisfies the
following estimate
\begin{equation}\label{5.34}
 \|u^g(t)\|_{\infty}\le  K_0+\int_0^t C_{qm} [\|\nabla g(\tau)\nabla u^g(\tau)\|_{q,loc}+
 \|\nabla g(\tau)\nabla u^g(\tau)\|_{m,loc}]d\tau.
\end{equation}

The proof of the  estimate can  easily be obtained by direct
computation from (\ref{5.33}) using the estimates of the Newton
potential given in lemma 5.1.

{\bf Lemma 5.5.} {\em Assume that conditions of lemma 5.4 hold.
Then given the function $u^g(t,x)$ of the form (\ref{5.33})  the
function $\nabla u^g(t,x)$ admits a representation of the form
$$
\nabla u^g(t,x)=E[\nabla u_0(\psi^g_{t,0}(x))\eta^{x,g}(t)-$$
\begin{equation}\label{5.35}
\int_0^t\frac 1{\sigma(t-\tau)}\nabla
p^g(\tau,\psi^g_{t,\tau}(x))\int_\tau^t\eta^{x,g}(\theta)d\hat
w(\theta) d\tau] \end{equation}
 and  the  estimate
$$\|\nabla u^g(t)\|_{\infty}\le e^{\int_0^t
K^1_g(\theta)d\theta}K_0^1+ $$\begin{equation}\label{5.36}
\int_0^tC_{qm}\frac 1{\sigma\sqrt{t-\tau}}e^{\int_\tau^t
K^1_g(\theta)d\theta}K^1_g(\tau)[\|\nabla
u^g(\tau)\|_{q,loc}+\|\nabla u^g(\tau)\|_{m,loc}]d\tau
\end{equation}
holds for $1<q<\frac 32<3<m<\infty$.}

Proof. To derive (\ref{5.35}) we compute directly the gradient of
the first term in (\ref{5.34}) and apply the Bismut-Elworthy -Li
formula \cite{EL} to compute the gradient of the second term in
this relation. To verify the estimate (\ref{5.36}) we use the
above estimates for the process $\eta(t)$  and the estimates of
the Newton potential derivative from lemma 5.1. Then we obtain
$$ \|\nabla u^g(t)\|_{\infty}\le e^{\int_0^t
K^1_g(\theta)d\theta}K_0^1+$$\begin{equation}\label{5.37} \int_0^t
C_{m,q}\frac 1{\sigma\sqrt{t-\tau}}e^{\int_\tau^t
K^1_g(\theta)d\theta}[\|\nabla g(\tau)\nabla
u^g(\tau)\|_{m,loc}+\end{equation}$$\|\nabla g(\tau)\nabla
u^g(\tau)\|_{q,loc}]d\tau]\le e^{\int_0^t
K^1_g(\theta)d\theta}K_0^1+$$$$\int_0^tC_{qm}\frac
1{\sigma\sqrt{t-\tau}}e^{\int_\tau^t
K^1_g(\theta)d\theta}K^1_g(\tau)[\|\nabla
u^g(\tau)\|_{q,loc}+\|\nabla u^g(\tau)\|_{m,loc}]d\tau].$$

Now we have to derive the estimate for  the function  $\|\nabla
u(t)\|_{r,loc}$.

{ \bf Lemma 5.6. } {\em Assume that the conditions of lemma 2.4
hold. Then for $1<r<\infty$ the function $ u^g(t,x)$  given by
(\ref{5.33}) satisfies the estimate
\begin{equation}\label{5.38}\|\nabla u^g(t)\|_{r,loc}\le   e^{\int_0^t
K^1_g(\theta)d\theta}\|\nabla u_0\|_{r,loc}+
\end{equation} $$C_{qm}\int_0^t e^{\int_0^\tau
K^1_g(\theta)d\theta}K^1_g(\tau)\|\nabla u^g (\tau)\|_{r,loc},
d\tau$$
 with a constant $C$ depending on $r$ and a certain
constant $T$ to be specified later.}

Proof.  Recall that  along with (\ref{5.35})  $\nabla u^g(t,x)$
admits the representation
$$\nabla u^g(t,x)=E[\nabla
u_0(\psi_{t,0}^g(x))\eta^{x,g}(t)-\int_0^t\nabla^2p^g(\tau,
\psi^g_{t,\tau}(x))\eta^{x,g}(\tau)d\tau].$$
To derive the estimate for $\|\nabla
u(t)\|_{r,loc}^r=\int_{G}\|\nabla u(t,x)\|^rdx$  (where $G$ is an
arbitrary  compact in $R^3$)  by the triangle inequality we get
$$\|\nabla u^g(t)\|_{r,loc}\le \alpha_1+\alpha_2,$$ where
$$\alpha_1=\left(\int_{G}E[\|\nabla
u_0(\psi^g_{t,0}(x))\eta^{x,g}(t)\|^r]dx\right)^{\frac 1r},$$
$$
\alpha_2=\left(\int_{G}\int_0^t\|\nabla^2
p^g(\tau,\psi^g_{t,\tau}(x))) \eta^{x,g}(\tau)\|^rd\tau
dx\right)^{\frac 1r}.$$

To estimate $\alpha_1$ we apply the H\"older inequality  and
recall that $\psi_{t,\tau}(x)$ preserves the volume. As a result
we have
$$
\alpha_1\le \left(\int_{G}(E[\|\nabla
u_0(\psi^g_{t,0}(x))\|^2]E[\|\eta^g(t)\|^2])^{\frac
r2}dx\right))^{\frac 1{r}}\le$$
$$ \|\nabla u_0\|_{r,loc} e^{\int_0^t
K^1_g(\theta)d\theta}.$$

To  estimate  $\alpha_2$ we apply  the
Calderon-Zygmund inequality and the above property of
$\psi_{t,\tau}(x)$ to obtain
$$\alpha_2^r\le C_r \int_0^t e^{\int_0^\tau K_g^1(\theta)d\theta}K^1_g(\tau)\int_{G}\|\nabla u(\tau,x)
\|^{r}dx d\tau.$$  Combining the above estimates for $\alpha_1$
and $\alpha_2$ we obtain the required  estimate

$$
\|\nabla u^g(t)\|_{r,loc}\le   e^{\int_0^t
K^1_g(\theta)d\theta}[\|\nabla u_0\|_{r,loc}+ $$$$C_r\int_0^t
e^{\int_0^\tau K_g^1(\theta)d\theta}K^1_g(\tau)\|\nabla
u^g(\tau)\|_{r,loc} d\tau]. $$

{\bf Theorem 5.7.} {\em Assume that conditions  {\bf C 5.1} and
{\bf C 5.2} hold.  Then there exists an interval
 $\Delta_1=[0,T_1]$ and functions $\alpha(t)$, $\beta(t)$, $\kappa$ bounded
for $t\in \Delta_1$, such that, if for all $t\in \Delta_1,$  $ \|g
(t)\|_{\infty}\le \kappa(t)$ and
 $ \|\nabla g(t)\|_{\infty
 }\le \alpha(t),$  $ \|\nabla g(t)\|_r\le \beta_r(t)$
 then the function $\|\nabla u^g(t,x)\|$ (where $u^g(t,x)$ is given by (\ref{5.21}))
 satisfies the estimates
\begin{equation}\label{5.39}
 \|u^g
(t)\|_{\infty}\le \kappa(t),\quad \|\nabla u^g(t)\|_{\infty
 }^2\le \alpha(t),\quad  \|\nabla u^g(t)\|_{r,loc}^2\le \beta_r(t)\end{equation}
for $r=q$ and $r=m$ and $1<m<\frac 32<3<q<\infty.$}

Proof. Analyzing the  above estimates for the functions $u^g(t,x)$
and $\nabla u^g(t,x)$ we get the following estimates
\begin{equation}\label{5.40}
\|\nabla u^g(t)\|_{\infty}\le e^{\int_0^t
K^1_g(\theta)d\theta}K_0^1+ \end{equation}
$$\int_0^t C_{qm}e^{\int_0^\tau K_g^1(\theta)d\theta}K^1_g(\tau)[\|\nabla
u^g(\tau)\|_{q,loc}+\|\nabla u^g(\tau)\|_{m,loc}]d\tau],$$

\begin{equation}\label{5.41}
\|\nabla u(t)\|_{r,loc}\le   e^{\int_0^t
K^1_g(\theta)d\theta}[\|\nabla u_0\|_{r,loc}+\end{equation}$$
C_{r}\int_0^t e^{\int_0^\tau
K_g^1(\theta)d\theta}K^1_g(\tau)\|\nabla u^g(\tau)\|_{r,loc}
 d\tau].$$

To derive the required estimates consider the integral equations

\begin{equation}\label{5.42}
\alpha(s)=e^{\int_s^t
\alpha(\theta)d\theta}K_0^1+C_{qm}\int_s^te^{\int_s^\tau
\alpha(\theta)d\theta} \alpha(\tau)[n_q(\tau)+n_m(\tau)]d\tau,
\end{equation}

$$n_{r}(s)=   e^{\int_s^t \alpha(\theta)d\theta}\|\nabla u_0\|_{r}+
C_{r}\int_s^te^{\int_s^\tau \alpha(\theta)d\theta}n_{r}(\tau)
\alpha(\tau)d\tau $$ for $r=q$ and $r=m$ and $C^1_{qm}=max(C_q,
C_m)$. Finally we consider the equation
\begin{equation}\label{5.43}
\beta(s)=   e^{\int_s^t \alpha(\theta)d\theta}\beta_0+
C^1_{qm}\int_s^te^{\int_s^\tau
\alpha(\theta)d\theta}\alpha(\tau)\beta(\tau)  d\tau,
\end{equation} where $\beta(\tau)=n_q(\tau)+n_m(\tau)$, and
$$\|\nabla u_0\|_{q,loc}+\|\nabla u_0\|_{m,loc}=n_q(0)+n_m(0)=\beta_0.$$
Next instead of the above system of integral equations we consider
the  system of ODEs
\begin{equation}\label{5.44}
\frac{d\alpha}{ds}=-\alpha^2(s)-{C_{qm}}\alpha(s)\beta(s),\quad
\alpha(t)=K_0^1, \end{equation}
\begin{equation}\label{5.45}
\frac{d\beta}{ds}=-\alpha(s)\beta (s)-
C^1_{qm}\alpha(s)\beta(s),\quad \beta (t)=\beta_0. \end{equation}
 By classical results of the ODE theory we know that there exists an interval
 $[0,T_1]$ depending on $K_0^1, N_0^1$ and
$C, C_{qm}$ such that the system (\ref{5.44}), (\ref{5.45}) has a
bounded solution defined on this interval.

To prove the convergence of   functions $u^k(t,x), \nabla
u^k(t,x)$  we need one more auxiliary estimate.  Actually, we have
proved that $u^k(t)\in Lip$ with the Lipschitz constant
independent of $k$. It remains to prove that $\nabla u^k(t)$ have
the same property.

 {\bf Lemma 5.8. } {\em Assume that {\bf C 5.1} and {\bf C 5.2} hold. Then
the function $\nabla u^g(t)$  satisfies the estimate
$$ \|\nabla u^g(t,x)-\nabla u^g(t,y)\|\le N_1^g(t)\|x-y\|^\alpha \, \mbox{  if  } t\in[0,T_1]$$
for any $x,y \in G$ where $G$ is a compact in $R^3$ and the
positive function $N_1^g(t)$ depending on parameters in conditions
{\bf C 5.1} and {\bf C 5.2}  is bounded over the interval $[0,T_1]
$ defined in theorem 5.7.}

Proof.  Applying the integration by parts Bismut -- Elworthy -- Li
formula to (\ref{5.33}) we deduce the following expression for the
gradient of the function $u(t,x)$
\begin{equation}\label{5.46}\nabla u^g(t,x)=E[\nabla
u_0(\psi^g_{t,0}(x))\eta^{x,g}(t)-\end{equation} $$\int_0^t\frac
1{\sigma(t-\tau)}\nabla
p^g(\tau,\psi^g_{t,\tau}(x))\int_\tau^t\eta^{x,g}(\theta)d \hat
w(\theta)d\tau].$$

It results from (\ref{5.46})  that
$$\|\nabla u^g(t,x)-\nabla u^g(t,y)\|\le
\kappa_1+\kappa_2+\kappa_3+\kappa_4,$$ where
$$\kappa_1= E[\|\nabla u_0(\psi^g_{t,0}(x))-\nabla u_0(\psi^g_{t,0}(y))\|\|\eta^{x,g}(t)\|],$$

$$\kappa_2= E[\|\nabla u_0(\psi^g_{t,0}(y))\|
 \|\eta^{x,g}(0)-\eta^{y,g}(0)\|],$$

$$
\kappa_3=\int_0^t\frac 1{\sigma(t-\tau)}E[\,\|\nabla p^g(\tau,
\psi^g_{t,\tau}(x))-$$$$\nabla p^g(\tau,
\psi^g_{t,\tau}(y))\|\|\int_\tau^t\eta^{x,g}(\theta)d \hat
w(\theta)\|\,]d\tau,$$
$$\kappa_4=\int_0^t\frac
1{\sigma(t-\tau)}E\left[\,\|\nabla p^g(\tau,
\psi^g_{t,\tau}(y))\|\int_\tau^t[\eta^{x,g}(\theta)-\eta^{y,g}(\theta)]d
\hat w(\theta)\|\,\right]d\tau.
$$

One can easily check using the estimates stated in lemmas 5.3 --
5.5 that  under conditions {\bf C 5.1}, {\bf C 5.2}
$$\kappa_1\le L^1_0
E\|\psi^g_{t,0}(x)-\psi_{t,0}(y)\|e^{\int_0^tK^1_g(\theta)d\theta}\le
\|x-y\|L^1_0e^{\int_\tau^tL_g(\theta)d\theta}$$ and
$$\kappa_2\le K_0^1E\|\eta^{x,g}(t)-\eta^{y,g}(t)\|\le \|x-y\|\int_\tau^tK^1_g(\theta)
 e^{\int_\theta^tK^1_g(\theta_1)d\theta_1}d\theta.$$

To derive the estimates for $\kappa_3$ and $\kappa_4$ we recall
(see lemma 5.1)  that the solution of the Poisson equation
$-\Delta p^g=\gamma^g$ satisfies the estimates
$\|\nabla_i\nabla_jp^g\|_\infty\le C
(\|\gamma^g\|_q+[\gamma^g]_{\alpha,G})$,
$\|\nabla_i\nabla_jp^g\|_r\le\|\gamma^g\|_r$  and $\|\nabla
p^g\|_\infty\le C_{qm}(\|\gamma^g\|_q+\|\gamma^g\|_m)$. Hence we
 obtain the inequalities
$$\kappa_3\le \int_0^t\frac 1{\sigma\sqrt{t-\tau}}
(E\|\psi^g_{t,\tau}(x)-\psi_{t,\tau}(y)\|^2)^{\frac
12}(\|\gamma^g(\tau)\|_q+[\gamma^g(\tau)]_{\alpha,
G})$$$$e^{\int_\tau^tK^1_g(\theta)d\theta}d\tau\le  \int_0^t\frac
1{\sigma\sqrt{t-\tau}} (\|x-y\|L^1_0e^{\int_\tau^t
L_g(\theta)d\theta}(\|\gamma^g(\tau)\|_q+$$
$$[\gamma^g(\tau)]_{\alpha,G})
e^{\int_\tau^tK^1_g(\theta)d\theta}d\tau$$ and
$$\kappa_4\le \int_0^t\frac
1{\sigma\sqrt{t-\tau}}C_{qm}(\|\gamma^g(\tau)\|_q+\|\gamma^g(\tau)\|_m)$$$$
(E\|\eta^{x,g}(\tau)-\eta^{y,g}(\tau)\|^2)^{\frac 12} d\tau\le$$
$$\|x-y\|\int_0^t\frac{C_{qm}(\|\gamma^g(\tau)\|_q+\|\gamma^g(\tau)\|_m)}{\sigma\sqrt{t-\tau}}
\int_\tau^tK^1_g(\theta)
 e^{\int_\theta^tK^1_g(\theta_1)d\theta_1}d\theta
 d\tau. $$
 Denote by $\Theta(t)=sup_{x,y\in G}\frac{\|\nabla u^g(t,x)-\nabla
 u^g(t,y)\|}{\|x-y\|^\alpha}$and note that
  $$[\gamma^g(\tau)]_{\alpha,G}=sup_{x,y\in G}[\frac{K^1_g(\tau)\|\nabla
u^g(\tau,x)-\nabla u^g(\tau,y)\|}{\|x-y\|^\alpha}+$$$$
\frac{K^1_u(\tau)\|\nabla g(\tau,x)-\nabla
g(\tau,y)\|}{\|x-y\|^\alpha}]=K^1_g(\tau)\Theta(\tau)+K^1_u(\tau)[g(\tau)]_{\alpha,G}.$$
Then combining the above estimates   for $\kappa_i,\, i=1,2,3,4,$
and applying the Gronwall lemma we derive the estimate
$$\Theta(t)\le N^g(t)e^{\int_0^tK^1_g(\tau)d\tau} =N_1^g(t),$$
where $N^g(t)$ is a positive bounded function defined on the
interval $[0,T_1]$ and depending on parameters in conditions {\bf
C 5.1} and {\bf C 5.2}.

The estimates of  theorem 5.7    and lemma 5.8 allow to prove the
uniform convergence on compacts of the successive approximations
(\ref{5.10})-(\ref{5.14}) for the solutions of the system
 (\ref{5.7}) -- (\ref{5.9})  in $C([0,T_1],$ \break $C^{1,\alpha}(K))\cap C([0,T_1], L^m(G)\cap L^q(G))$
 for  $1<q<\frac 32<3<m<\infty$   and arbitrary compact $G$ in $R^3$.

To this end we differentiate the system (\ref{5.10})-(\ref{5.14})
and add to this system the following relations
\begin{equation}\label{5.47}
d\eta^{k,x}_{t,\theta}=-\nabla
u^k(\theta,\psi^k_{t,\theta})\eta^{x,k}_{t,\theta}d\theta,\quad
\eta^{x,k}_{t,t}=I,
 \end{equation}
where $I$ is the identity matrix acting in  $R^3$ and
$$
\nabla u^{k+1}(t,x)= E[\nabla
u_0(\psi^{k+1}_{t,0}(x))\eta^{x,k}_{t,0}-$$\begin{equation}\label{5.48}\int_0^t\frac
1{\sigma(t-\tau)}\nabla
p^{k+1}(\tau,\psi^k_{t,\tau}(x))\int_\tau^t\eta_{t,\theta}^{x,k}d\hat
w(\theta) d\tau], \end{equation}
\begin{equation}\label{5.49}
-2 \nabla p^{k+1}(t,x)=\int_0^\infty \frac 1\tau E[
\gamma^{k+1}(t,x+B(\tau))B(\tau)]d\tau,
 \end{equation}
 where $\gamma^{k+1}=\nabla u^{k+1}\nabla u^k.$

Now we can prove the following assertion.

 {\bf Theorem 5.9. } {\em Assume that  {\bf C 5.2} holds.
 Then  if $k\to\infty$ the functions $u^k(t), \nabla u_k(t,x)$ determined by  (\ref{5.8}) and
 (\ref{5.48})
  uniformly
converge on compacts to a limiting  function $u(t)\in C([0,T_1],
C^{1,\alpha}), 0<\alpha\le 1$ for all
 $t\in [0, T_1]$, where   $[0,T_1]$  is the interval such that the solution of  (\ref{5.45}), (\ref{5.46})
 is bounded on $[0,T_1]$. In addition on
 this interval  the limiting function satisfies the
 estimates
 $\sup_x\|\nabla u(t,x)\| \le \alpha(t)$ , $\|\nabla u(t)\|_{q, loc} \le \beta(t)$
 for $1<q<\frac 32$ where $\alpha(t)$ and $\beta(t)$ solve (\ref{5.45}), (\ref{5.46}).}

Proof. By theorem 5.7 we know that the mapping
 $$\Phi(t,x,g)=E\left[u_0(\psi^g_{t,0}(x))-\int_0^t\nabla p^g(\tau,\psi^g_{t,\tau}(x))d\tau\right]$$ acts in the space
 $  C^{1,\alpha}\cap L_{q,loc}\cap L_{m,loc} $ (for a fixed $t\in [0.T_1])$  with  $1<q<\frac 32<3<m<\infty.$

Consider the successive approximations (\ref{5.10}) --(\ref{5.14})
and (\ref{5.47}) -- (\ref{5.49}),  denote by
$$S^{k+1}(t,x)=\|u^{k+1}(t,x) -u^k(t,x)\|,$$
$$ n^{k+1}(t,x)= \|\nabla u^{k+1}(t,x) -\nabla u^k(t,x)\|
$$
and let
$$
l^k(t)=\|S^k(t)\|_\infty,\quad m_r^k(t)= \|S^k(t)\|_r,$$
$$ \rho^k(t)=\|n^k(t)\|_\infty,  \quad \zeta^k_r(t)=\|n^k(t)\|_r.$$

Then we obtain
$$n^{k+1}(t,x)\le L_0^1  (E[\|\psi^k_{t,0}(x)-\psi^{k-1}_{t,0}(x)\|\|\eta^{x,k}_{t,0}\|]
+$$$$E[\|\psi^k_{t,0}(x)\|\|\eta^{x,k}_{t,0}-\eta^{x,k-1}_{t,0}\|])
+\int_0^t\frac 1{\sigma (t-\tau)}E[\|\nabla p^{k+1}(\tau,
\psi^k_{t,\tau}(x))-$$$$\nabla p^{k}(\tau,
\psi^{k-1}_{t,\tau}(x))\|\|\int_\tau^t\eta^{x,k}_{t,\theta}d\hat
w(\theta)\|]d\tau+$$
\begin{equation}\label{5.50}\int_0^t\frac 1{\sigma (t-\tau)}E\left[\|\nabla p^{k}(\tau,
\psi^k_{t,\tau}(x))\|\int_\tau^t[\eta^{x,k}_{t,\theta}-\eta^{x,k-1}_{t,\theta}]d\hat
w(\theta)\|\right]d\tau.\end{equation}

Recall that by lemmas 5.2, 5.3  we know that
$$\sup_xE\|\psi^k_{t,0}(x)-\psi^{k-1}_{t,0}(x)\|\le
\int_0^t[\,\|u^k(\tau)- u^{k-1}(\tau)\|_\infty ]d\tau
e^{\int_0^t\alpha(\tau)d\tau},$$
$$\sup_x E\|\eta^{x,k}_{t,0}-\eta^{x,k-1}_{t,0}\|\le
\int_0^t\|\nabla u^k(\tau)-\nabla u^{k-1}(\tau)\|_\infty d\tau
e^{\int_0^t\alpha(\tau)d\tau}
$$
$$
+\sup_x\int_0^t E\|\nabla u^{k-1}(\tau,\psi^k_{t,\tau}(x))-\nabla
u^{k-1}(\tau,\psi^{k-1}_{t,\tau}(x))\| d\tau
e^{\int_0^t\alpha(\tau)d\tau}$$ and applying the estimates from
theorem 5.7 we get
$$\rho^{k+1}(t)\le
e^{\int_0^t\alpha(\tau)d\tau}[L_0^1\int_0^t\sup_xE\|u^k(\tau,\psi^k_{t,\tau}(x))
-u^{k-1}(\tau,\psi^{k-1}_{t,\tau}(x))\| d\tau$$$$
+\int_0^t\rho^k(\tau)d\tau +\sup_x\int_0^t E\|\nabla
u^{k-1}(\tau,\psi^k_{t,\tau}(x))-\nabla
u^{k-1}(\tau,\psi^{k-1}_{t,\tau}(x))\| d\tau ]
$$$$+\int_0^t\frac 1{\sigma \sqrt{t-\tau}}C[\|\nabla u^k(\tau)\nabla u^{k-1}(\tau)\|_q+
\|\nabla u^k(\tau)\nabla u^{k-1}(\tau)\|_m]$$$$
(E\|\eta^k(\tau)-\eta^{k-1}(\tau)\|^2_\infty)^{\frac 12} d\tau+$$
$$
\int_0^t\frac {e^{\int_\tau^t
\alpha(\theta)d\theta}}{\sigma \sqrt{t-\tau}}\sup_xE\|\nabla p^{k+1}(\tau,
\psi^k_{t,\tau}(x))-\nabla p^{k}(\tau,
\psi^{k-1}_{t,\tau}(x))\|^2)^{\frac 12}d\tau.$$

To derive the estimate for the last term  we recall ( see lemma
5.1) that for $1<q<\frac 32$  the inequality
 $$\|\nabla
p^k(t,x)-\nabla p^k(t,y)\|\le
\|\nabla^2p^k(t)\|_\infty\|x-y\|\le$$$$ C[\|\gamma
^k(t)\|_{q,loc}+[\gamma^k(t)]_{\alpha,G}]\|x-y\|$$ holds and as a
result we obtain

$$ E\|\nabla p^{k}(\tau,
\psi^k_{t,\tau}(x))-\nabla p^{k}(\tau,
\psi^{k-1}_{t,\tau}(x))\|\le$$$$ C[\beta(\tau)+\Theta(\tau)]E\|
\psi^k_{t,\tau}(x)- \psi^{k-1}_{t,\tau}(x)\|. $$

In addition
$$\|\nabla
p^{k+1}(t)-\nabla p^k(t)\|_\infty\le
C_{qm}[\,\|\gamma^{k+1}(t)-\gamma^k(t)\|_{q,loc}
+$$$$\|\gamma^{k+1}(t)-\gamma^k(t)\|_{m,loc}]\le C_{qm}
\alpha(t)[\,\|\nabla u^{k+1}(t)-\nabla
u^{k}(t)\|_{q,loc}+$$$$\|\nabla u^{k}(t)-\nabla
u^{k-1}(t)\|_{q,loc}+$$$$\|\nabla u^{k+1}(t)-\nabla
u^{k}(t)\|_{m,loc}+\|\nabla u^{k}(t)-\nabla
u^{k-1}(t)\|_{m,loc}].$$

It results from  (\ref{5.50}) that
$$n^{k+1}(t,x)\le
C(t)[\int_0^tE\|\nabla u^k(\tau,\psi^k_{t,\tau}(x)) -\nabla
u^{k-1}(\tau,\psi^{k-1}_{t,\tau}(x))\|
d\tau+$$$$\int_0^tn^k(\tau,x)d\tau] +\int_0^t\frac 1{\sigma
\sqrt{t-\tau}}C_1[\|\nabla u^k(\tau)\nabla u^{k-1}(\tau)\|_q+$$$$
\|\nabla u^k(\tau)\nabla u^{k-1}(\tau)\|_m]^r
(E\|\eta^{x,k}(\tau)-\eta^{x,k-1}(\tau)\|^{2})^{\frac 12} d\tau$$
$$
+\int_0^t\frac 1{\sigma \sqrt{t-\tau}} e^{\int_\tau^t
\alpha(\theta)d\theta}(E\|\nabla p^{k+1}(\tau,
\psi^k_{t,\tau}(x))-\nabla p^{k}(\tau,
\psi^{k-1}_{t,\tau}(x))\|^2)^{\frac 12} d\tau.$$

Note that by the H\"older inequality we can prove that for any
positive $f(\tau)\in L^r $ and $\frac 1{m_1}+\frac 1r=1$
$$\int_G[\int_0^tf(\tau,x)d\tau]^rdx\le \int_G t^{\frac
r{m_1}}\int_0^tf^r(\tau,x)d\tau dx=$$$$t^{\frac
r{m_1}}\int_0^t\int_Gf^r(\tau,x)dx d\tau$$ and    for $ \frac
1{m_1}+\frac 1r=1$ and $m_1<2$ we have
\begin{equation}\label{5.51}\int_G[\int_0^t\frac
1{\sigma\sqrt{t-\tau}}f(\tau,x)d\tau]^rdx\le
t^{\frac{r(2-m_1)}{2m_1}}
\int_0^t\int_Gf^r(\tau,x)dxd\tau.\end{equation}

Then from  (\ref{5.50})  and (\ref{5.51})  we  have for $r>2$
$$\zeta^{k+1}_r(t)\le
C_2(t)[\int_0^t\int_G [E\|u^k(\tau,\psi^k_{t,\tau}(x))
-u^{k-1}(\tau,\psi^{k-1}_{t,\tau}(x))\|^rdx d\tau]
+$$$$\int_0^t\zeta^k_{r}(\tau)d\tau+\int_0^t \int_G\|\nabla
u^{k-1}(\tau,\psi^k_{t,\tau}(x))-\nabla
u^{k-1}(\tau,\psi^{k-1}_{t,\tau}(x))\|^r dx d\tau]
$$$$+\int_0^t\frac 1{\sigma \sqrt{t-\tau}}C[[\|\nabla u^k(\tau)\nabla u^{k-1}(\tau)\|_q+
\|\nabla u^k(\tau)\nabla u^{k-1}(\tau)\|_m]^r$$$$
\int_G(E\|\eta^{x,k}(\tau)-\eta^{x,k-1}(\tau)\|^{2})^{\frac r2}
dx] d\tau +\int_0^t\frac 1{\sigma \sqrt{t-\tau}} e^{\int_\tau^t
\alpha(\theta)d\theta}$$$$\int_G (E\|\nabla p^{k+1}(\tau,
\psi^k_{t,\tau}(x))-\nabla p^{k}(\tau,
\psi^{k-1}_{t,\tau}(x))\|^2)^{\frac r2}dx d\tau.$$

In addition  for $m^k_r(t)=\|u^k(t)-u^{k-1}(t)\|_{r,loc}$  using the
apriori estimates proved in lemmas 5.2 -- 5.8 and theorem  5.9 we
obtain

$$m^{k+1}_r(t)\le C(t)[(\int_0^t\int_G E\|u^k(\tau,\psi^k_{t,\tau}(x))-u^{k-1}(\tau,
\psi^{k-1}_{t,\tau}(x))\|^rdxd\tau)^{\frac 1r}
$$
$$+(\frac 1\sigma t^{\frac 1{m_1}-\frac 12}\int_0^t \int_G
E\|\nabla u^{k+1}(\tau, \psi^{k+1}_{t,\tau}(x))\nabla u^{k}(\tau,
\psi^{k}_{t,\tau}(x))-$$$$\nabla u^{k}(\tau,
\psi^{k}_{t,\tau}(x))\nabla u^{k-1}(\tau,
\psi^{k-1}_{t,\tau}(x))\|^rdx d\tau)^{\frac 1r}]\le$$
$$ C_1(t)[(\int_0^tm^k_r(\tau)d\tau)^{\frac 1r}+(\int_0^t\int_G
\alpha(\tau)E\|\psi^k_{t,\tau}(x)-\psi^{k-1}_{t,\tau}(x)\|^rdxd\tau)^{\frac
1r}+$$$$\frac 1\sigma t^{\frac 1{m_1}-\frac 12}(\int_0^t
[\rho^{k+1}(\tau)+\rho^{k}(\tau)]\zeta_r^{k}(\tau)d\tau)^{\frac
1r}].$$

Since   $u^k$ and $\nabla u^k$ are uniformly bounded on $[0,T_1]$
and
 $$\|\nabla u^1(t,\cdot)-\nabla u_0(\cdot)\|_{r,loc}\le const <\infty,$$
both for $r=m$ and $r=q$ we obtain that there exists a bounded on
$[0,T_1]$ positive function $C_2(t)$ such that  the function
$\kappa^n(t)=\rho^n(t)+\zeta^n_m(t)+m^n_r$  satisfies  the
estimate
$$\kappa^n(t)\le \frac{[C_2(t)]^n}{n!}$$ and hence
$\lim_{n\to\infty}\kappa^n(t)=0,$ since $C_2(t)$ is bounded on
$[0,T_1]$. Finally we obtain that
 for each  $t\in [0,T_1)$
  the family $u^n(t,\cdot)$  uniformly converges
to a limiting function  $u(t,\cdot)\in C^{1,\alpha}\cap L_{m,
loc}$. In addition,
 we can check that the limiting function  $\nabla u(t,x)$ is Lipschitz continuous in $x$.
In fact, by lemma 2.8 and theorem 2.9 for each $t\in [0,T_1]$,
 we have for any $x,y\in G$
  $$\|\nabla u^n(t,x)-\nabla u^n(t,y)\|\le N(t)\|x-y\|,$$
where  $N(t)$  and $T_1$ were defined above in lemmas 5.8 and
theorem 5.7 respectively  and the estimate is uniform in $n$. This
allows to state that
  the limiting function is Lipschitz continuous as well.

To prove  the uniqueness of the solution of (2.8)-(2.10)
constructed above we assume first that there exist two solutions
 $u_1(t,x)$, $u_2(t,x)$  to
(\ref{5.7})-(\ref{5.9}) possessing the same initial data
    $u_1(0,x)=u_2(0,x)= u_0(x)$.

Computations similar to those used to prove the convergence of the
family $(u^n(t),\nabla u^n(t))$ allow to check that
$$[\nabla u_1(t)-\nabla u_2(t)]_{\alpha, G}=0 \quad \mbox{
and } \quad \|\nabla u_1(t)-\nabla u_2(t)\|_{m,loc}=0.$$

Finally, we know that the Cauchy problemfor a stochastic equation
 with Lipschitz coefficients has a  unique solution.
This implies the uniqueness  of the solution to
(\ref{5.7})-(\ref{5.9}).

Summarizing the above results we see that the following statement
is valid.

{\bf Theorem 5.10.}  {\em Assume that  {\bf C 5.2} holds. Then
there exists a unique solution $\psi_{t,x}(s), u(t,x), p(t,x)$ of
the system (\ref{5.7})-(\ref{5.9}), for all $t$ from the interval
the $[0,T_1]$, with $T_1$ given by theorem 5.7 and $x\in G$ for
any compact $G\subset R^3$. In addition 
$\psi_{t,x}(s)$ is a Markov process in $R^3$ and $u\in C([0,T_1],
C^{1,\alpha}(G))\cap C([0,T_1], L_{q,loc}\cap L_{m,loc})$ for
$1<q<\frac 32<3<m<\infty.$}

To fulfill our program we have to check that the conditions of
theorem 2.8 are sufficient to  verify  that the functions
$u(t,x),p(t,x)$ given by  (\ref{5.8}), (\ref{5.9}) define a weak
solution of the Navier -Stokes system.

To this end we have to apply  the results of  the Kunita theory of
stochastic flows.

Namely we  check that given a distribution valued process
$\lambda(t)$ of the form
\begin{equation}\label{5.52}\lambda(t)=u_0-\int_0^t\nabla p^u(\tau)\circ\phi^u_{0,\tau}d\tau
\end{equation}

 the function
$$\lambda(t)\circ\psi^u_{t,0}= u_0\circ\psi^u_{t,0}-\int_0^t\nabla p^u(\tau)\circ
 \psi^u_{t,\tau}d\tau$$
gives rise to a solution of (\ref{5.1}).

To this end we apply the generalized Ito formula \cite{Ku1},
\cite{Ku2} to derive

\begin{equation}\label{5.53}
\lambda(t)\circ\psi^u_{t,0}=u_0+\int_{0}^{t} \frac
{\sigma^2}2\Delta[u(\theta)\circ \psi^u_{\theta,0}]d\theta+
\end{equation}
$$\int_{0}^{t} \nabla[
u(\theta)\circ\psi^u_{\theta,0}] \sigma dw(\theta)- \int_{0}^{t}
\nabla[ u(\theta)\circ\psi^u_{\theta,0}] u(\theta)
d\theta-\int_0^t\nabla p^u(\theta)d\theta.
$$

Note that for $Lu=-(u,\nabla)u+ \frac {\sigma^2}2\Delta u$ we have

 $$E\left[\int_{R^3}\int_0^t(L(u(\tau)\circ\psi^u_{\tau, s}(x))d\tau,
 h(x))dx\right]=$$
$$ E\left[\int_0^t\langle u(\tau)\circ\psi^u_{\tau,0}, L^*h\rangle
 d\tau\right]=\int_0^tL\langle E[u(\tau\circ\psi^u_{\tau,0})], h\rangle
 d\tau.$$

Hence
$$
u(t)= E[\lambda(t)\circ \psi^u_{t,0}]= u_0+\int_0^tLE[u(\tau)\circ
\psi^u_{\tau,0}]d\tau-\int_0^t\nabla p^u(\tau)d\tau.$$

Differentiating each term with respect to $t$ we can check that
the function
\begin{equation}\label{5.54}
u(t)=E[\lambda(t)\circ\psi^u_{t,0}] \end{equation}
 solves the
Cauchy problem (\ref{5.1}),(\ref{5.2}).

To summarize the obtained results  we can state the following
assertion.

{\bf Theorem 5.11. }{\em Assume that {\bf C 5.2} holds. Then the
functions $u(t,x),$ $ p(t,x)$ given by (\ref{5.8}),(\ref{5.9}) are
defined on the interval $[0,T_1]$ with $T_1$ determined by theorem
5.8 and satisfy (\ref{5.1})-(\ref{5.2}) in a weak sense on this
interval.}

{\bf Remark 5.12.}  We have proved that under condition {\bf C
5.2} the system (\ref{5.7})-({\ref{5.9}) gives rise to a weak
solution of (\ref{5.1})-(\ref{5.2}) . Note that if the initial
data are  smoother, say $u_0\in C^{2+\alpha}$, $\alpha\in [0,1]$
similar considerations can be applied to verify that the pair
$u(t,x), p(t,x)$ given by (\ref{5.8})-(\ref{5.9})  stands for a
classical $C^2$-smooth solution of (\ref{5.1})-(\ref{5.2}).  In
fact in this case applying the generalized Ito formula for the
verification assertion we may treat the action of the operator
$ L$ in the classical sense rather then in the weak sense.

\section{Lagrangian and stochastic approach to the Euler and the N-S
system}

\setcounter{equation}{0}
 The probabilistic approach  developed in
the previous section is in a sense an analogue of the Lagrangian
approach to the Euler and the Navier-Stokes systems. A rather
close model was constructed in papers by Constantin and Iyer
\cite{Con1},\cite{CI1}. To make it easier to compare we rewrite
the results from these papers in terms similar to those used in
the previous section.  We consider first
 the Euler system
\begin{equation}\label{6.1}
 \frac{\partial u}{\partial t}+(u,\nabla )u=-\nabla p ,\quad u(0,x)=u_0(x),\quad x\in R^3
 \end{equation}
\begin{equation}\label{6.2}
 div\,u=0.\end{equation}
and recall that  the corresponding Lagrangian path starting at $y$
is governed by the Newton equation
\begin{equation}\label{6.3}\frac{\partial^2 \tilde\phi_{0,t}(y)}{\partial
t^2}=F_{\tilde\phi}(t,y).\end{equation}

The incompressibility condition for the map $\phi$ yields
\begin{equation}\label{6.4}{ \det}(\nabla
\tilde\phi_{0,t}(y))=1.\end{equation}

The force $F$ in (3.3) has the form
\begin{equation}\label{6.5}
F_{\tilde\phi}(t,y)=-\nabla p(t,\tilde\phi_{0,t}(y))= -[(\nabla
\tilde\phi_{0,t}(y))^*]^{-1}\nabla [p(t,\tilde\phi_{0,t}(y))].
\end{equation} One can deduce from (\ref{6.3}) that
\begin{equation}\label{6.6}
\frac {\partial}{\partial t}[\frac {\partial\tilde\phi^k_{0,t}(y)}
{\partial t}\frac{\partial\tilde\phi^k_{0,t}(y)}{\partial
y_i}]=-\frac {\partial  q(t,\tilde\phi_{0,t}(y))}{\partial
y_i},\end{equation}
 where
\begin{equation}\label{6.7}q(t,y)= p(t,y)-\frac 12 \| \frac
{\partial\tilde\phi_{0,t}(y)}{\partial t}\|^2. \end{equation} We recall that in
(\ref{6.6}) and below   summation over the repeated
indices is assumed. Integrating (\ref{6.6}) in time we get
\begin{equation}\label{6.8}
\frac {\partial\tilde\phi^k_{0,t}(y)}{\partial t}\frac
{\partial\tilde\phi^k_{0,t}(y)}{\partial y_i}=u_0(y)-
\frac{\partial
 n(t,\tilde\phi_{0,t}(y))}{\partial y_i},
\end{equation}
 where
\begin{equation}\label{6.9}
 n(t,y)=\int_0^t q (\tau,y)d\tau
\end{equation}
 and
\begin{equation}\label{6.10}
u_0(y)=\frac {\partial\tilde\phi_{0,t}(y)}{\partial t}|_{t=0}
 \end{equation}
 is the initial velocity.

 Consider the inverse diffeomorphism $\tilde\psi_{t,0}=[\tilde\phi_{0,t}]^{-1}$,
 come back to (\ref{6.7}), multiply it by
$[\nabla\tilde\psi_{t,0}]$ and put $y=\tilde\psi_{t,0}(x)$. As a
result we obtain by the chain rule the relation
\begin{equation}\label{6.11}
u^i(t,x)=(u_0^j(\tilde\psi_{t,0}(x))\nabla_{x_i}\tilde\psi^j_{t,0}(x)-\int_0^t\nabla_{x_i}
q(\tau,\tilde\psi_{t,\tau}(x))d\tau. \end{equation}

The equation (\ref{6.11})  shows that the general Euler velocity
may be written in the form that generalizes the Clebsch variable
representation
$$u=[\nabla \tilde\psi_{t,0}]^*C-\nabla n,$$
where $C=u_0(\psi_{t,0}(x))$ is an active vector and $n$ is
defined by the incompressibility condition $div u=0$.

Note  that a vector $A$ is called active if
$$ \frac d{dt}A=\frac{\partial A}{\partial t}+(u,\nabla)A=0.$$
 It is easy to check by the chain rule that
\begin{equation}\label{6.12}\frac d{dt}\tilde\psi_{t,\theta}(x)=\frac {\partial\tilde\psi_{t,\theta}(x)}{\partial
 t}+(u,\nabla)\tilde\psi_{t,\theta}(x)=0,\end{equation}
 that is $\psi_{t,0}(x)$ is an active vector.

 Hence the Euler equations are equivalent to the system consisting of
 (3.9) and the following relation

\begin{equation}\label{6.13}
\Delta n(t,x)= \frac{\partial }{\partial
x_i}\{u_0^k(\tilde\psi_{t,0}(x))\frac{\partial\tilde\psi^k_{t,0}(x)
}{\partial
 x_i}\},
\end{equation}
 where $n$ is given by (\ref{6.9}).

Now one can assume the periodic boundary conditions  or   the zero
boundary conditions at infinity.
 Note that in the periodic case
$ n(t,x ), u(t,x)$ and
\begin{equation}\label{6.14}
\delta(t,x)=x-\tilde\psi_{t,0}(x)
 \end{equation}
 are periodic functions in each spatial direction.  Finally  due to
 $div u=0$ one can rewrite  the
 equation of state  (\ref{6.11})   in the
 form
 \begin{equation}\label{6.15}
u(t)={\Pi}\{u_0^j(\tilde\psi_{t,0})\nabla\tilde\psi^j_{t,0}\}=
{\Pi}\{[\nabla\tilde\psi_{t,0}]^*u_0(\tilde\psi_{t,0})\},
 \end{equation}
where $ {\Pi}=I-\nabla\Delta^{-1}\nabla$ is the Leray-Hodge
projector (with corresponding boundary conditions) on divergence
free vector fields. The Euler pressure is determined up to
additive constants by
$$
p(t,x)=\frac{\partial n(t,x)}{\partial t}+ (u(t,x),\nabla)
n(t,x)+\frac 12 \|u(t,x)\|^2.
$$
Note that (\ref{6.11}), (\ref{6.12}) made a closed system and may
be used to determine $u(t)$.

Let us  compare (\ref{6.11}), (\ref{6.12}) with the alternative
representation for the state $u(t)$ developed in the previous
section.

To this end we choose  $\phi_{0,t}:y\to \phi_{0,t}(y)$ to be  a
volume preserving diffeomorphism  that satisfies the equation
\begin{equation}\label{6.16}d\phi_{0,\tau}(y)=u(t-\tau,\phi_{0,\tau}(y))d\tau,
\quad \phi_{0,0}(y)=y, \end{equation} with div $u(t)=0$.

Consider the system
\begin{equation}\label{6.17}d\psi_{t,\theta}(x)=-u(\theta,
 \psi_{t,\theta}(x)) d\theta, \quad \psi_{t,t}(x)=x,
  \end{equation}
\begin{equation}\label{6.18}
u(t,x)=u_0(\psi_{t,0}(x))-\int_0^t\nabla p(\tau,
\psi_{t,\tau}(x))d\tau,  \end{equation}
\begin{equation}\label{6.19}
-2p(t,x)=E[\int_0^\infty\gamma(t, x+B(\tau))d\tau],
\end{equation} where $\gamma$  is given by (\ref{1.4}).

If the fields $u(t,x), p(t,x)$ are regular enough then we may
construct the representation of the solution to the Euler system
in the form  (\ref{6.18}), (\ref{6.19}).

To check this we consider a volume preserving diffeomorphism
$\phi_{0,t}:y\to \tilde\phi_{0,t}(y)$ that satisfies (\ref{6.16}).

Next we consider the system
\begin{equation}\label{6.20} d\psi_{t,\theta}(x)=-u(\theta,
 \psi_{t,\theta}(x)) d\theta, \quad \psi_{t,t}(x)=x,
 \end{equation}
\begin{equation}\label{6.21}
u(t,x)=u_0(\psi_{t,0}(x))-\int_0^t\nabla p(\tau,
\psi_{t,\tau}(x))d\tau, \end{equation}
\begin{equation}\label{6.22}
-2p(t,x)=E[\int_0^\infty\gamma(t, x+B(\tau))d\tau], \end{equation}
where $\gamma$  is given by (\ref{1.3}).

If $u(t,x)$ is regular enough then we may construct the
representation of the solution to the Euler system (\ref{6.1}),
(\ref{6.2})
 in the form (\ref{6.21}), (\ref{6.22}).

To this end we consider a vector field $\lambda(t)$  satisfying
the equation
$$ \frac{d\lambda(t)}{dt}= -\nabla p(t)\circ\phi_{0,t}\quad
\lambda(0)=u_0, $$ where $\phi_{0,t}$ satisfies the ODE
(\ref{6.16}), and let the process $\psi_{t,0}$ be its inverse.
Applying the Kunita approach \cite{Ku1} to the  process
$\psi_{t,0}$ we can verify that $\psi_{t,0}$ along with
(\ref{6.20}) satisfies the equation
\begin{equation}\label{6.23}
\psi_{t,\tau}(x)=x+\int_\tau^t\nabla
 \phi^g_{\theta,t}(\psi_{t,\theta})^{-1}u(\theta,x)d\theta,
 \end{equation}
that allows to prove that $u(t)$ given by (\ref{6.21}) satisfies
(\ref{6.1}).

 Comparing (\ref{6.15}) and (\ref{6.21}) we note that
 they give different expressions for the velocity field. Actually
 (\ref{6.21}) includes the Euler pressure  $p(t,x)$ instead of $q(t,x)$ used in
 (\ref{6.15}).
  Besides  the probabilistic representation of the
  solution $p$ to the Poisson equation
  $$-\Delta p=\nabla_iu_k\nabla_ku_i$$ is used instead of
  the Leray projection.

Coming back to the Navier-Stokes system ((\ref{5.1}),(\ref{5.2})
we recall here the approach due to Constantin and Iyer \cite{CI1}.
The stochastic counterpart of the Navier-Stokes equations in the
version of Iyer \cite{I} looks like the following.

 Consider the closed
stochastic system
\begin{equation}\label{6.24}d\phi_{0,\theta}=u(\theta, \phi_{0,\theta})dt+\sigma
dw(\theta),\quad \phi_{0,0}(y)=y,
\end{equation}
\begin{equation}\label{6.25} \psi_{\theta,0}=[\phi_{0,\theta}]^{-1},
 \end{equation}
\begin{equation}\label{6.26}
u(t)=E {\Pi}[(\nabla \psi_{t,0})(u_0\circ\psi_{t,0})].
 \end{equation}
The existence and uniqueness of the solution to this system is
proved in \cite{CI1} by the successive approximation technique. As
a result the authors  constructed  a strong local in time solution
of the Cauchy problem for the Navier-Stokes system for regular
enough initial data.

The main result due to Constantin and Iyer reads as follows

{\bf Theorem 6.1.} {\em  Let $k\ge 1$ and $u_0\in C^{k+1,\alpha}$
be divergence free. Then there exists a time interval $[0,T]$ with
$T=T(k,\alpha, L, \|u_0\|_{k+1,\alpha})$ but independent of
viscosity $\sigma$ and a pair $\phi_{0,t}(x), u(t,x)$  such that
$u\in C([0,T], C^{k+1,\alpha})$ and $(u,\phi)$ satisfy
(\ref{6.24})-(\ref{6.26}). Further there exists  $U=U(k,\alpha,
L,\|u_0\|_{k+1,\alpha})$  such that $\|u(t)\|_{k+1,\alpha}\le U$
for $t\in [0,T]$  and $u$ satisfies the N-S system.}

As we have mentioned above an  approach close to the one of \cite{CI1} was developed in our
previous paper \cite{AB}.  Both these approaches allow to construct
a classical (local in time) solution to the Cauchy problem for the
Navier-Stokes system and prove the uniqueness of the solution.

On the other hand the approach developed in section 5 allows to
construct a weak (local in time)  solution to (\ref{5.1}),
({\ref{5.2}) and prove the uniqueness of this solution in the
corresponding functional classes.

The  stochastic counterpart of the Navier-Stokes system considered
in section 2 has the form

\begin{equation}\label{6.27} d\psi_{t,\theta}(x)=-u(\theta,
 \psi_{t,\theta}(x) d\theta +\sigma d\hat w(\theta), \quad \psi_{t,t}(x)=x,
  \end{equation}
\begin{equation}\label{6.28}
u(t,x)= E[u_0(\psi_{t,0}(x))-\int_0^t\nabla
p(\tau,\psi_{t,\tau}(x))d\tau], \end{equation}

\begin{equation}\label{6.29}2p(t,x)=-\int_0^\infty E[\gamma(t,
x+B(\tau))]d\tau.  \end{equation}

Note that we can use the relation
\begin{equation}\label{6.30}
-2\nabla p(t,x)=E[\int_0^\infty\frac 1\tau\gamma(t,
x+B(\tau))B(\tau)d\tau]
\end{equation}
to eliminate the pressure from the above system  (\ref{6.27}) --
(\ref{6.29}).

We can see that the difference between (\ref{6.27}) --
(\ref{6.29}) and (\ref{6.24}) -- (\ref{6.26}) has the same nature
as the difference between (\ref{6.11}), (\ref{6.12}) and
(\ref{6.17}) -- (\ref{6.19}).

Finally we note that the approach developed in section 5 allows us
to construct both strong (classical) and  weak (distributional)
solutions of the Cauchy problem for the N-S system.

{\bf Acknowledgement.} The authors gratefully acknowledge the
financial support of DFG Grant 436 RUS 113/823.

\end{document}